# Improvements on Waring’s Problem

Li An-Ping

Beijing 100085, P.R. China
apli0001@sina.com

Abstract

By a new recursive algorithm for the auxiliary equation, in this paper, we will give some improvements for Waring’s problem.

## 1. Introduction

Waring's problem is now to find $G(k)$, the least integer $s$, such that each sufficient large integer may be represented a sum of at most $s$ $k$th powers of natural numbers. The Hardy-Littlewood method, that is, so called circle method is main analysis method, which is proposed by Hardy. Ramanujian and Littlewood in about 1920's, which have been applied successfully in solving some problems of number theory, e.g. Waring's problem and Goldbach' s problem.
The known best results for Waring's problem up to now are as following

For sufficiently large $k$ (Wooley [5] ),

$$G(k) \le k\left(\log(k\log k) + O(1)\right). \tag{1.1}$$

And for smaller $k$ ,

$$G(5) \le 17, G(6) \le 24, G(7) \le 33, G(8) \le 42, \cdots. \tag{1.2}$$

For the details is referred to see the Vaughan and Wooley's survey paper [4]. In this paper, by a new recursive algorithm, we will give some improvements for $G(k)$.

**Theorem 1.1.** For sufficiently large $k$,

$$G(k) \le 4k + o(k). \tag{1.3}$$

**Theorem 1.2.** For $5 \le k \le 20$, let $F(k)$ be as in the List 1.1, then

$$G(k) \le F(k). \tag{1.4}$$

| $k$ | $F(k)$ | $k$ | $F(k)$ | $k$ | $F(k)$ | $k$ | $F(k)$ |
|---|---|---|---|---|---|---|---|
| 5 | 17 | 9 | 33 | 13 | 49 | 17 | 66 |
| 6 | 21 | 10 | 37 | 14 | 51 | 18 | 68 |
| 7 | 25 | 11 | 41 | 15 | 57 | 19 | 75 |
| 8 | 32 | 12 | 45 | 16 | 64 | 20 | 75 |

List 1.1

In section 4, there is further progress.

## 2. The Proof of Theorems 1.1.

Suppose that $P$ is a sufficient large integer, $\mathscr{C}(P)$ is a subset of $[0,P]$, $k$ is a given integer, consider the equation

$$x_1^k+\cdots+x_s^k=y_1^k+\cdots+y_s^k,\quad x_i,y_i\in\mathscr{C}(P),1\le i\le s. \tag{2.1}$$

Denoted by $S_s(\mathscr{C}(P))$ the number of solutions of (2.1), or simply $S_s(P)$, when the selection $\mathscr{C}(P)$ is clear in context. The equation above is called auxiliary equation of Waring problem.

In the following, we will take use of iterative method to construct $\mathscr{C}(P)$.

Suppose that $\theta$ is a real number, $0<\theta<1$. Let $\tilde{P}=P^{1+\theta}$, $\mathscr{P}$ is a set of prime numbers $p$ in interval $[P^\theta/2,P^\theta]$, write $|\mathscr{P}|=Z$, define

$$\mathscr{C}(\tilde{P})=\left\{x\cdot p\mid x\in\mathscr{C}(P),p\in\mathscr{P}\right\}. \tag{2.2}$$

With respect to the construction, we will also consider following a relative equation

$$\begin{aligned}&p^k\left(x_1^k+\cdots+x_{s-i}^k-y_1^k+\cdots+y_{s-i}^k\right)=q^k\left(x_{s-i+1}^k+\cdots+x_s^k-y_{s-i+1}^k+\cdots+y_s^k\right),\\&x_i,y_i\in\mathscr{C}(P),1\le i\le s.\ \text{where}\ p,q\in\mathscr{P},p\ne q.\end{aligned} \tag{2.3}$$

Denote by $T_{s,i}(p,q)$ the number of solutions of (2.3), and $T_{s,i}(q)=\sum_p T_{s,i}(p,q)$.

**Lemma 1**. For integer $i, 0<i<s$, it has

$$S_s(\tilde{P})\ll Z^{2s-1}T_{s,i}(q). \tag{2.4}$$

Proof. As usual, write $e(x)=e^{2\pi ix}$, let

$$f(\alpha)=\sum_{x\in\mathscr{C}(P)}e\left(x^k\alpha\right),\ f(\alpha,p)=\sum_{x\in\mathscr{C}(P)}e\left(p^kx^k\alpha\right),\ \tilde{f}(\alpha)=\sum_{y\in\mathscr{C}(\tilde{P})}e\left(y^k\alpha\right).$$

Then clearly,

$$\tilde{f}(\alpha)=\sum_{p\in\mathscr{P}}\sum_{x\in\mathscr{C}(P)}e\left(p^kx^k\alpha\right)=\sum_{p\in\mathscr{P}}f(\alpha,p).$$

Applying Hölder's inequality, it has

$$S_s(\tilde{P})=\int_0^1\left|\tilde{f}(\alpha)^{2s}\right|d\alpha=\int_0^1\left|\tilde{f}(\alpha)^s\,\tilde{f}(\alpha)^{-s}\right|d\alpha$$

$$=\int_0^1\left(\sum_{p\in\mathscr{P}}\sum_{q\in\mathscr{P}}f(\alpha,p)f(-\alpha,q)\right)^s d\alpha$$

$$= \int_0^1 \left( \sum_{p=q} + \sum_{p\neq q} \right)^s d\alpha$$

$$\ll \int_0^1 \left( Z\cdot |f(\alpha)|^2 \right)^s d\alpha + \int_0^1 ( \sum_{p,q\in\mathscr{P},\ p\neq q} |f(\alpha,p)f(-\alpha,q)| )^s d\alpha$$

$$\ll Z^s S_s(P) + Z^{2(s-1)} \sum_{p,q\in\mathscr{P},\ p\neq q} \int_0^1 |f(\alpha,p)f(\alpha,q)|^s d\alpha.$$

It is clear that $Z^s S_s(P)$ is minor for $S_s(\tilde{P})$.

Moreover, for a non-negative integer $i$, let $\Lambda_i(\alpha,p,q) = f(\alpha,p)^{s-i} f(\alpha,q)$, then by Cauchy inequality, it has

$$\int_0^1 \left| f(\alpha,p)^s f(\alpha,q)^s \right| d\alpha = \int_0^1 |\Lambda_i(\alpha,p,q)||\Lambda_i(\alpha,q,p)| d\alpha$$

$$\le \left( \int_0^1 |\Lambda_i(\alpha,p,q)|^2 d\alpha \right)^{1/2} \left( \int_0^1 |\Lambda_i(\alpha,q,p)|^2 d\alpha \right)^{1/2}$$

And

$$\sum_{p,q\in\mathscr{P},p\neq q} \int_0^1 \left| f(\alpha,p)^s f(\alpha,q)^s \right| d\alpha \le \sum_{p,q\in\mathscr{P},p\neq q} \left( \int_0^1 |\Lambda_i(\alpha,p,q)|^2 d\alpha \right)^{1/2} \left( \int_0^1 |\Lambda_i(\alpha,q,p)|^2 d\alpha \right)^{1/2}$$

$$\le \sum_{p,q\in\mathscr{P},p\neq q} \int_0^1 |\Lambda_i(\alpha,p,q)|^2 d\alpha$$

Clearly, the inner integral is the number of solutions of equation (2.3).

□

Denote by $\mathscr{P}[a,b]$ the set of prime numbers in the interval $[a,b]$. Similar to construct (2.2), let $\theta_1,\cdots,\theta_k$ be $k$ real numbers, $0\le\theta_i\le 1/k,\ 1\le i\le k,$ which will be determined later, and let $Z_i = (P/2)^{\theta_i}$, $\mathscr{P}_i = \mathscr{P}[Z_i/2, Z_i]$, $P_{i+1} = P_i / Z_{i+1}$, $H_i = P/Z_i^k$, $1\le i\le k$. Recursively define

$$\mathscr{C}(P_i) = \left\{ x\cdot p \mid x\in\mathscr{C}(P_{i+1}),\ p\in\mathscr{P}_{i+1}, (p,x)=1 \right\}, \quad i=0,1,\cdots k.$$

Simply write $\mathscr{C}_i = \mathscr{C}(P_i), \quad i=0,1,\cdots,k.$

In the following, it will be used the notation of difference of a function: As usual, for an integer coefficient of polynomial $\phi(x)$, recursively define the forward differences

$$\Delta_1(\phi(x),t) = \phi(x+t) - \phi(x),$$

$$\Delta_{i+1}(\phi(x),h_1,h_2,\cdots,h_i,h_{i+1}) = \Delta_1\left(\Delta_i(\phi(x),h_1,h_2,\cdots,h_i),h_{i+1}\right). \quad i=1,2,\cdots.$$

Suppose that $t = h\cdot m$, $m$ is a constant, then we know that $m \mid (\phi(x+t)-\phi(x))$, in this case we define modified differences

$$\Delta_1^*(\phi(x),h;m)=m^{-1}\left(\phi(x+hm)-\phi(x)\right).$$

$$\Delta_{i+1}^*(\phi(x),h_1,\cdots,h_i,h_{i+1};m_1,\cdots,m_i,m_{i+1})=\Delta_1^*\left(\Delta_i^*(\phi(x),h_1,\cdots,h_i;m_1,\cdots,m_i),h_{i+1};m_{i+1}\right).$$

Simply write,

$$\Psi_i=\Psi_i(x,h_1,\cdots,h_i;p_1^k,\cdots,p_i^k)=\Delta_i^*(x^k,h_1,\cdots,h_i;p_1^k,\cdots,p_i^k),\qquad i=1,2,\cdots.$$

And define

$$f_i(\alpha)=\sum_{x\in\mathscr{C}_i}e(x^k\alpha),\quad f_i(\alpha,p)=\sum_{x\in\mathscr{C}_i,(x,p)=1}e(p^kx^k\alpha),\quad p\in\mathscr{P}_i,\qquad 1\le i\le k.$$

$$g_i(\alpha;h_1,\cdots,h_i;p_1^k,\cdots,p_i^k)=\sum_x e(\Psi_i(x)\alpha),\qquad 1\le i\le k. \tag{2.5}$$

And

$$F_i(\alpha,q)=\sum_{1\le h_1\le H_1}\cdots\sum_{1\le h_i\le H_i}\sum_{p_1\in\mathscr{P}_1}\cdots\sum_{p_i\in\mathscr{P}_i}g_i(\alpha q^k)$$

$$F_i^{(2)}(\alpha,q)=\sum_{1\le h_1\le H_1}\cdots\sum_{1\le h_i\le H_i}\sum_{p_1\in\mathscr{P}_1}\cdots\sum_{p_i\in\mathscr{P}_i}\left|g_i(\alpha q^k)\right|^2 \tag{2.6}$$

$$F_i^{(4)}(\alpha,q)=\sum_{1\le h_1\le H_1}\cdots\sum_{1\le h_i\le H_i}\sum_{p_1\in\mathscr{P}_1}\cdots\sum_{p_i\in\mathscr{P}_i}\left|g_i(\alpha q^k)\right|^4$$

Let

$$\mathcal{L}_0=\sum_p\int_0^1|f(\alpha,p)|^{2(s-2)}|f(\alpha,q)|^4\,d\alpha$$

$$\mathcal{L}_i=\int_0^1|f_i(\alpha)|^{2(s-2)}F_i^{(2)}(\alpha,q)d\alpha,\qquad 1\le i<k. \tag{2.7}$$

$$\mathcal{J}_i=\int_0^1|f_i(\alpha)|^{2(s-2)}F_i(\alpha,q)d\alpha,\quad 1\le i<k.$$

**Lemma 2.2**

$$\int_0^1\sum_{p\in\mathscr{C}_{i+1}}|f_{i+1}(\alpha,p)|^{2(s-2)}F_i^{(4)}(\alpha,q)d\alpha$$

$$\ll(Z_{i+1})^{k+1}P^2\tilde{Z}_i\tilde{H}_iS_{s-2}(P_{i+1})+P(Z_{i+1})^k\mathcal{J}_{i+1}+P\mathcal{L}_{i+1}. \tag{2.8}$$

Proof. As usual, for a number $x$, denote by $\overline{x}\in\mathbb{Z}_{p^k}$ with that $x\equiv\overline{x}\mod p^k$. Simply write $\Psi(x)=\Psi_i(x)$, and $\sum_{\boldsymbol{p}}=\sum_{p_1\in\mathscr{P}_1}\cdots\sum_{p_i\in\mathscr{P}_i}$, $\sum_{\boldsymbol{h}}=\sum_{h_1}\cdots\sum_{h_i}$.

For a $p\in\mathscr{P}_{i+1}$, let

$$\mathscr{D}_i(\Psi,p)=\{x\mid x\in\mathscr{C}_i,\Psi'(x)\equiv 0\mod p\},\quad \Omega(p)=\mathscr{C}_i^2\setminus\mathscr{D}_i(\Psi,p)^2.$$

Then the sum of $g_i(\alpha q^k)^2$ can be divided two parts, that is, normal and singular parts, i.e. one with $(x,y)\in\Omega(p)$, and the other one not. It is not difficult to demonstrate that the integral in the singular part is secondary, for the simplicity, we save the investigation, and in the following acquiesce in the normal part.

It is clear that for equation $\Psi(\overline{x})+\Psi(\overline{y})\equiv n \mod p^k$ there are at most $O(p^k)$ solutions with that $(x,y)\in\Omega(p)$, hence we can divide $\Omega(p)$ into $O(p^k)$, say $l$ classes, $\{\mathscr{F}_i(p)\}_1^l$, such that in each class $\mathscr{F}_i(p)$ the equation has at most two solutions $\mod p^k$. And denote by $\beta_j(x)=card\{y\mid(x,y)\in\mathscr{F}_j(p)\}$. It is clear that $\sum_{1\le j\le l}\beta_j(x)\le P.$ Hence, it has

$$\int_0^1\sum_{p\in\mathscr{C}_{i+1}}|f_{i+1}(\alpha,p)|^{2(s-2)}F_i^{(4)}(\alpha,q)d\alpha=\sum_{\boldsymbol{p}}\sum_{\boldsymbol{h}}\int_0^1\sum_{p\in\mathscr{C}_{i+1}}|f_{i+1}(\alpha,p)|^{2(s-2)}|g_i(\alpha,q)|^4\,d\alpha$$

$$\ll\sum_{\boldsymbol{p}}\sum_{\boldsymbol{h}}\sum_{p\in\mathscr{C}_{i+1}}\int_0^1|f_{i+1}(\alpha,p)|^{2(s-2)}\left|\sum_{1\le j\le l}\sum_{(\overline{x},\overline{y})\in\mathscr{F}_j(p)}e((\Psi_i(x)+\Psi_i(y))q^k\alpha)\right|^2d\alpha$$

$$\ll l\cdot\sum_{\boldsymbol{p}}\sum_{\boldsymbol{h}}\sum_{p\in\mathscr{C}_{i+1}}\sum_{1\le j\le l}\int_0^1|f_{i+1}(\alpha,p)|^{2(s-2)}\left|\sum_{(\overline{x},\overline{y})\in\mathscr{F}_j(p)}e((\Psi_i(x)+\Psi_i(y))q^k\alpha)\right|^2d\alpha$$

$$\ll l\cdot\sum_{\boldsymbol{p}}\sum_{\boldsymbol{h}}\sum_{p\in\mathscr{C}_{i+1}}\sum_{1\le j\le l}\int_0^1|f_{i+1}(\alpha,p)|^{2(s-2)}\times\Bigg(|\mathscr{F}_j|+\sum_{h}\sum_{x}\beta_j(x)e(\Delta_1^*(\Psi_i(x),h;p)q^k\alpha)$$

$$+\sum_{(\overline{x},\overline{y})\in\mathscr{F}_j}\sum_{h_1,h_2}e(\Delta_1^*(\Psi_i(x),h_1;p)+\Delta_1^*(\Psi_i(y),h_2;p))q^k\alpha)\Bigg)d\alpha$$

$$\ll l\sum_{\boldsymbol{p}}\sum_{\boldsymbol{h}}\sum_{p\in\mathscr{C}_{i+1}}\sum_{1\le j\le l}\int_0^1|f_{i+1}(\alpha,p)|^{2(s-2)}|\mathscr{F}_j(p)|d\alpha$$

$$+l\sum_{\boldsymbol{p}}\sum_{\boldsymbol{h}}\sum_{p\in\mathscr{C}_{i+1}}\int_0^1|f_{i+1}(\alpha,p)|^{2(s-2)}\sum_{1\le j\le l}\sum_{h}\sum_{x}\beta_j(x)e(\Delta_1^*(\Psi_i(x),h;p))q^k\alpha)d\alpha$$

$$+l\sum_{\boldsymbol{p}}\sum_{\boldsymbol{h}}\sum_{p\in\mathscr{C}_{i+1}}\int_0^1|f_{i+1}(\alpha,p)|^{2(s-2)}\sum_{1\le j\le l}\sum_{(\overline{x},\overline{y})\in\mathscr{F}_j}\sum_{h_1,h_2}e(\Delta_1^*(\Psi_i(x),h_1;p)+\Delta_1^*(\Psi_i(y),h_2;p))q^k\alpha)d\alpha$$

$$\ll\sum_{\boldsymbol{p}}\sum_{\boldsymbol{h}}lP^2Z_{i+1}\int_0^1|f_{i+1}(\alpha)|^{2(s-2)}\,d\alpha$$

$$+l\cdot\sum_{\boldsymbol{p}}\sum_{\boldsymbol{h}}\int_0^1|f_{i+1}(\alpha,p)|^{2(s-2)}P\sum_{p\in\mathscr{C}_{i+1}}\sum_{h}\sum_{x}e(\Delta_1^*(\Psi_i(x),h;p))q^k\alpha)d\alpha$$

$$+l\sum_{\boldsymbol{p}}\sum_{\boldsymbol{h}}\int_0^1|f_{i+1}(\alpha)|^{2(s-2)}\sum_{p\in\mathscr{C}_{i+1}}\left|\sum_{h}\sum_{x}e(\Delta_1^*(\Psi_i(x),h;p)q^k\alpha)\right|^2d\alpha$$

$$\ll lP^2Z_{i+1}\tilde{Z}_i\tilde{H}_iS_{s-2}(P_{i+1})$$

$$+P\sum_{1\le j\le l}\int_0^1|f_{i+1}(\alpha)|^{2(s-2)}\sum_{\boldsymbol{p}}\sum_{\boldsymbol{h}}\sum_{p\in\mathscr{C}_{i+1}}\sum_{h}\sum_{x}e(\Delta_1^*(\Psi_i(x),h;p))q^k\alpha)d\alpha$$

$$+lH_{i+1}\int_0^1|f_{i+1}(\alpha,p)|^{2(s-2)}\sum_{\boldsymbol{p}}\sum_{\boldsymbol{h}}\sum_{p\in\mathscr{C}_{i+1}}\sum_{h}\left|\sum_{x}e(\Delta_1^*(\Psi_i(x),h;p)q^k\alpha)\right|^2d\alpha$$

$$\ll lP^2 Z_{i+1}\tilde{Z}_i\tilde{H}_i S_{s-2}(P_{i+1}) + Pl\mathcal{J}_{i+1} + P\mathcal{L}_{i+1}$$

□

From the proof above, we can know that

$$\mathcal{L}_0 = T_{s,2}(q) = \sum_p \int_0^1 |f(\alpha,p)|^{2(s-2)} |f(\alpha,q)|^4 d\alpha$$
$$\ll Z^k P^2 Z S_{s-2}(P) + Z^k P\mathcal{J}_1(P) + P\mathcal{L}_1 \tag{2.9}$$

Besides,

$$\mathcal{J}_i = \int_0^1 |f_i(\alpha)|^{2(s-2)} F_i(\alpha,q)d\alpha \le \left(\int_0^1 |f_i(\alpha)|^{2(s-2)} d\alpha\right)^{1/2} \left(\int_0^1 |f_i(\alpha)|^{2(s-2)} F_i(\alpha,q)^2 d\alpha\right)^{1/2}$$
$$\le S_{s-2}(P_i)^{1/2}(\tilde{H}_i\tilde{Z}_i)^{1/2}\mathcal{L}_i^{1/2}. \tag{2.10}$$

where $\tilde{H}_i = \prod_{j\le i} H_j$, $\tilde{Z}_i = \prod_{j\le i} Z_j$ .

In general, we have

**Lemma 2.3.**

$$\mathcal{L}_i \ll U_i + V_i + W_i. \tag{2.11}$$
$$U_i = \left(S_{s-2}(P_i)\right)^{1/2} Z_{i+1}^{(2s-5)/2} \times (\tilde{Z}_i\tilde{H}_i)^{1/2} \times (lP^2 Z_{i+1}\tilde{Z}_i\tilde{H}_i S_{s-2}(P_{i+1}))^{1/2},$$
$$V_i = \left(S_{s-2}(P_i)\right)^{1/2} Z_{i+1}^{(2s-5)/2} \times (\tilde{Z}_i\tilde{H}_i)^{1/2} \times (Pl\mathcal{J}_{i+1})^{1/2},$$
$$W_i = \left(S_{s-2}(P_i)\right)^{1/2} Z_{i+1}^{(2s-5)/2} \times (\tilde{Z}_i\tilde{H}_i)^{1/2} \times (P\mathcal{L}_{i+1})^{1/2}.$$

Proof.

$$\mathcal{L}_i = \int_0^1 |f_i(\alpha)|^{2(s-2)} F_i^{(2)}(\alpha,q)d\alpha \le \left(\int_0^1 |f_i(\alpha)|^{2(s-2)} d\alpha\right)^{1/2} \left(\int_0^1 |f_i(\alpha)|^{2(s-2)} (F_i^{(2)}(\alpha,q))^2 d\alpha\right)^{1/2}$$
$$\le \left(S_{s-2}(P_i)\right)^{1/2} \left(\int_0^1 \left|\sum_{p\in\mathscr{P}_{i+1}} f_{i+1}(\alpha,p)\right|^{2(s-2)} (F_i^{(2)}(\alpha,q))^2 d\alpha\right)^{1/2}$$
$$\le \left(S_{s-2}(P_i)\right)^{1/2} \left(Z_{i+1}^{2s-5}\int_0^1 \sum_{p\in\mathscr{P}_{i+1}} |f_{i+1}(\alpha,p)|^{2(s-2)} (F_i^{(2)}(\alpha,q))^2 d\alpha\right)^{1/2}$$
$$\le \left(S_{s-2}(P_i)\right)^{1/2} Z_{i+1}^{(2s-5)/2} \times \left(\tilde{Z}_i\tilde{H}_i\int_0^1 \sum_{p\in\mathscr{P}_{i+1}} |f_{i+1}(\alpha,p)|^{2(s-2)} F_i^{(4)}(\alpha,q)d\alpha\right)^{1/2}$$
$$\le \left(S_{s-2}(P_i)\right)^{1/2} Z_{i+1}^{(2s-5)/2} \times (\tilde{Z}_i\tilde{H}_i)^{1/2} \times \left((lP^2 Z_{i+1}\tilde{Z}_i\tilde{H}_i S_{s-2}(P_{i+1})) + (Pl\mathcal{J}_{i+1}) + (P\mathcal{L}_{i+1})\right)^{1/2}$$
$$\ll U_i + V_i + W_i .$$

□

Let

$$U_{i-1} = V_{i-1} \cdot H_i^{-\alpha/4} P^{-\tau/4} + W_{i-1} \cdot H_i^{-\alpha/2} P^{-\tau/2}, \quad \alpha < 2\,, \ \alpha + \tau < 2\,. \tag{2.12}$$

$\alpha, \tau$ will be decided later . Hence, it has

$$\left((Z_i)^k P^2 Z_i \tilde{Z}_{i-1} \tilde{H}_{i-1} S_{s-2}(P_i)\right)^{1/2} =$$

$$\left(P(Z_i)^k S_{s-2}(P_i)^{1/2} (\tilde{H}_i \tilde{Z}_i)^{1/2} \mathcal{L}_i^{1/2}\right)^{1/2} (H_i)^{-\alpha/4} P^{-\tau/4} + (P\mathcal{L}_i)^{1/2} (H_i)^{-\alpha/2} P^{-\tau/2}$$

And it follows

$$\mathcal{L}_i^{1/4} = \frac{\left(\left((Z_i)^k P^2 Z_i \tilde{Z}_{i-1} \tilde{H}_{i-1} S_{s-2}(P_i)\right)^{1/2} P^{1/2} (H_i)^{-\alpha/2} P^{-\tau/2}\right)^{1/2}}{P^{1/2} (H_i)^{-\alpha/2} P^{-\tau/2}}.$$

That is,

$$\mathcal{L}_i = (Z_i)^k Z_i \tilde{Z}_{i-1} \tilde{H}_{i-1} S_{s-2}(P_i) P (H_i)^{\alpha} P^{\tau}. \tag{2.13}$$

On the other hand,

$$\begin{aligned} \mathcal{L}_i &= U_i + V_i + W_i \approx U_i \cdot (H_{i+1}^{\alpha})^{1/2} P^{\tau/2} \\ &= \left(S_{s-2}(P_i)\right)^{1/2} Z_{i+1}^{(2s-5)/2} \times (\tilde{Z}_i \tilde{H}_i)^{1/2} \times (lP^2 Z_{i+1} \tilde{Z}_i \tilde{H}_i S_{s-2}(P_{i+1}))^{1/2} (H_{i+1})^{\alpha/2} P^{\tau/2}. \end{aligned} \tag{2.14}$$

Combine the two equalities,

$$\begin{aligned} &(Z_i)^k Z_i \tilde{Z}_{i-1} \tilde{H}_{i-1} S_{s-2}(P_i) P (H_i)^{\alpha} P^{\tau} \\ &= \left(S_{s-2}(P_i)\right)^{1/2} Z_{i+1}^{(2s-5)/2} \times (\tilde{Z}_i \tilde{H}_i)^{1/2} \times (lP^2 Z_{i+1} \tilde{Z}_i \tilde{H}_i S_{s-2}(P_{i+1}))^{1/2} (H_{i+1})^{\alpha/2} P^{\tau/2}. \end{aligned}$$

It follows

$$\frac{\left(S_{s-2}(P_i)\right)^{1/2}}{\left(S_{s-2}(P_{i+1})\right)^{1/2}} = \frac{Z_{i+1}^{(2s-4)/2} \times (H_i) \times (l)^{1/2} (H_{i+1})^{\alpha/2}}{\left((Z_i)^k\right) (H_i)^{\alpha} P^{\tau/2}}.$$

i.e.

$$\frac{S_{s-2}(P_i)}{S_{s-2}(P_{i+1})} = \frac{Z_{i+1}^{(2s-4)} \times (H_i)^2 \times (Z_{i+1})^k (H_{i+1})^{\alpha}}{(Z_i)^{2k} (H_i)^{2\alpha} P^{\tau}}. \tag{2.15}$$

Let $S_{s-2}(X) = X^{\lambda_{s-2}}$, (2.15) becomes

$$\frac{P_i^{\lambda_{s-2}}}{P_{i+1}^{\lambda_{s-2}}} = Z_{i+1}^{(2s-4)} \frac{(P/Z_i^k)^{2-2\alpha} (P/Z_{i+1}^k)^{\alpha} (Z_{i+1})^k}{(Z_i)^{2k} P^{\tau}}.$$

i.e.

$$(Z_{i+1})^{\lambda_{s-2}} = Z_{i+1}^{(2s-4)} P^{2-\alpha-\tau} Z_{i+1}^{(1-\alpha)k} Z_i^{-k(4-2\alpha)}.$$

And

$$\theta_i = \frac{\left((2s-4) + (1-\alpha)k - \lambda_{s-2}\right)}{k(4-2\alpha)} \theta_{i+1} + \frac{1}{2k} \frac{2-\alpha-\tau}{2-\alpha}.$$

Denote by $a = \dfrac{\left((2s-4)+(1-\alpha)k-\lambda_{s-2}\right)}{k(4-2\alpha)}$, $b = \dfrac{1}{2k}\dfrac{2-\alpha-\tau}{2-\alpha}$, it has

$\theta_i = a^j\theta_{i+j} + b\dfrac{1-a^j}{1-a}$, and $\theta_{k-i} = a^i\theta_k + b\dfrac{1-a^i}{1-a}$.

Besides, by (2.12), it has

$$\left((Z_k)^k P^2 Z_k \tilde{Z}_{k-1}\tilde{H}_{k-1}S_{s-2}(P_k)\right)^{1/2} = \left(P\mathcal{L}_k\right)^{1/2}(H_k)^{-\alpha/2}P^{-\tau/2}$$

It is easy to know that $\mathcal{L}_k = \tilde{Z}_k\tilde{H}_kP^2S_{s-2}(P_k)$, so it has

$$(Z_k)^k P^2 Z_k\tilde{Z}_{k-1}\tilde{H}_{k-1}S_{s-2}(P_k) = P\tilde{Z}_k\tilde{H}_kP^2S_{s-2}(P_k)(H_k)^{-\alpha}P^{-\tau}.$$

And $(H_k)^{2-\alpha}P^{-\tau} = 1 \Rightarrow (P/(Z_k)^k)^{2-\alpha}P^{-\tau} = 1 \Rightarrow \theta_k = \dfrac{2-\alpha-\tau}{k(2-\alpha)}$.

Define $\Delta(s) = \lambda_s + k - 2s$, there is

$$\theta_{k-i} = a^i\theta_k + b\frac{1-a^i}{1-a} = \frac{2-\alpha-\tau}{(2-\alpha)k+\Delta(s-2)} + \left(\frac{2-\alpha-\tau}{k(2-\alpha)} - \frac{2-\alpha-\tau}{(2-\alpha)k+\Delta(s-2)}\right)a^i.$$

Especially,

$$\theta = \theta_1 = \frac{2-\alpha-\tau}{(2-\alpha)k+\Delta(s-2)} + \left(\frac{2-\alpha-\tau}{k(2-\alpha)} - \frac{2-\alpha-\tau}{(2-\alpha)k+\Delta(s-2)}\right)a^{k-1}. \qquad (2.16)$$

On the other hand, by Lemma 2.1 with $i = 2$, it has

$$S_s(\tilde{P}) \ll Z^{2s-1}\mathcal{L}_0(P) \approx Z^{2s-1}U_0H^{\alpha/2}P^{\tau/2} = Z^{2s-1}Z^kP^2ZS_{s-2}(P)H^{\alpha/2}P^{\tau/2}.$$

and

$$\lambda_s = \frac{\lambda_{s-2}+2+(\alpha+\tau)/2}{(1+\theta)} + \frac{\theta(2s+k-k\alpha/2)}{(1+\theta)}.$$

Or,

$$\Delta(s) = \frac{\Delta(s-2)}{(1+\theta)} + \frac{\theta k-1}{(1+\theta)}(2-\alpha/2) + \frac{\tau}{2(1+\theta)}. \qquad (2.17)$$

When $k$ is greater, it may has

$$\theta \doteq \frac{(2-\alpha-\tau)}{(2-\alpha)k+\Delta(s-2)}. \qquad (2.18)$$

And let

$$(2-\alpha) = \frac{\beta\cdot\Delta(s-2)}{k}, \quad (2-\alpha-\tau) = \frac{\omega\cdot\Delta(s-2)}{k}. \qquad (2.19)$$

Substituting (2.18) and (2.19) in (2.17), it follows

$$\Delta(s)=\Delta(s-2)\left(1-\frac{3\omega}{2(k(1+\beta)+\omega)}\right)-\frac{k\rho}{k(1+\beta)+\omega},$$

where $\rho=(1+\beta-\omega)$. And then

$$\Delta(s)=\Delta(s-2i)\left(1-\frac{3\omega}{2(k(1+\beta)+\omega)}\right)^i-\frac{2\rho k}{3\omega}\left(1-\left(1-\frac{3\omega}{2(k(1+\beta)+\omega)}\right)^i\right).$$

So,

$$\Delta(d+2i)=\left(\Delta(d)+\frac{2\rho k}{3\omega}\right)\left(1-\frac{3\omega}{2(k(1+\beta)+\omega)}\right)^i-\frac{2\rho k}{3\omega}. \qquad (2.20)$$

$d=1$, or $2$, or others. It is known that $\lambda_1=1, \lambda_2=2$, i.e. $\Delta(1)=k-1, \Delta(2)=k-2$. Hence,

$$\Delta(\varepsilon+2i)=\left(k-\varepsilon+\frac{2\rho k}{3\omega}\right)\left(1-\frac{3\omega}{2(k(1+\beta)+\omega)}\right)^i-\frac{2\rho k}{3\omega}. \qquad \varepsilon=1, or 2.$$

Let

$$\Im(\beta,\omega)=\log\left(\frac{2\rho k}{3\omega(k-\varepsilon)+2\rho k}\right)\Big/\log\left(1-\frac{3\omega}{2(k(1+\beta)+\omega)}\right). \qquad (2.21)$$

From (2.20), we can know that $\Delta(2i+\varepsilon)$ will approach zero as $i$ tends to $\Im(\beta,\omega)$. Hence it has

**Lemma 2.4.** For sufficient large $k$, and arbitrary small $\varepsilon>0$, there is $s\le 2(\Im(\beta,\omega)+1)$ such that $\Delta(s)<\varepsilon.$

Besides, for greater $k$, it has

$$\Im(\beta,\omega)\sim\frac{2(k(1+\beta)+\omega)}{3\omega}\log\left(1+\frac{3\omega}{2\rho}\right)\to k \quad \text{as} \quad \omega\to 0. \qquad (2.22)$$

Moreover, we know that (see [3] or [1])

$$G(k)\le 3+2u+2\left[\frac{\Delta(u)}{2\hat{\sigma}}\right].$$

where $\hat{\sigma}=\dfrac{\log(1+1/k)}{4(1+\lambda)}$, $\lambda=\log k+\log\log k+O\left(\dfrac{\log\log k}{\log k}\right)$.

Take $u=2(\Im(\beta,\omega)+1)$, and let $\omega\to 0$, Theorem 1.1 is proved.

## 3. The proof of Theorem 1.2.

For the smaller $k$, $\lambda_s$ and $\Delta(s)$ may be followed by recursion (2.16) and (2.17) from the initial $\lambda_d$, $d=1,2,$ or any known better $\lambda_d$ by choosing optimal values $\alpha$ and $\tau$ in turn. With involve searching of two parameters $\alpha$ and $\tau$ , and the restriction of the ability of PC, we have lessen the search range only in four digits. The results in List 3.1 are obtained by PC.

| $k$ | $s$ | $\Delta(s)$ | $k$ | $s$ | $\Delta(s)$ |
|---|---|---|---|---|---|
| 5 | 8 | 0.000042 | 13 | 24 | 0.000014 |
| 6 | 10 | 0.000014 | 14 | 25 | 0.000016 |
| 7 | 12 | 0.000058 | 15 | 28 | 0.000009 |
| 8 | 14 | 0.000000 | 16 | 28 | 0.000002 |
| 9 | 16 | 0.000017 | 17 | 32 | 0.000019 |
| 10 | 18 | 0.000000 | 18 | 33 | 0.000012 |
| 11 | 20 | 0.000016 | 19 | 36 | 0.000008 |
| 12 | 22 | 0.000000 | 20 | 36 | 0.000005 |

List 3.1

For $k>2$, let $\varsigma(k)=4k$ if $k$ is a power of $2$, or $3k/2$ else. From the known results (see [2], [3]), we know that for two positive integers $t,v$, if satisfying

i) $2t+v\geq\varsigma(k),$

ii) $v\cdot 2^{1-k}>\Delta(t),$

then

$$G(k)\leq 2t+v. \tag{3.2}$$

With List 3.1, we take $v$ as in the following list

| $k$ | 5 | 6 | 7 | 8 | 9 | 10 | 11 | 12 | 13 | 14 | 15 | 16 | 17 | 18 | 19 | 20 |
|---|---|---|---|---|---|---|---|---|---|---|---|---|---|---|---|---|
| $v(k)$ | 1 | 1 | 1 | 4 | 1 | 1 | 1 | 1 | 1 | 1 | 1 | 8 | 2 | 2 | 3 | 3 |

List 3.2

And Theorem 1.2 is followed.

## 4. Further Improvements

Shortly after the paper appeared, we realize that the method of parameterized recursion applied in the sections 2,3 is also available for the recursive process applied in paper [1], and it is unexpected that the results are even better than the previous ones, the new results are that

**Theorem 4.1.** For sufficient large $k$ ,

$$G(k) \le \begin{cases} 3.661k + o(k), & \text{if } k \text{ is not a power of 2.} \\ 4\text{k}, & \text{else.} \end{cases} \tag{4.1}$$

**Theorem 4.2.** For $5 \le k \le 20$, let $F(k)$ be as in the List 4.1, then

$$G(k) \le F(k). \tag{4.2}$$

| $k$ | $F(k)$ | $k$ | $F(k)$ | $k$ | $F(k)$ | $k$ | $F(k)$ |
|---|---|---|---|---|---|---|---|
| 5 | 17 | 9 | 29 | 13 | 43 | 17 | 55 |
| 6 | 21 | 10 | 33 | 14 | 45 | 18 | 57 |
| 7 | 23 | 11 | 37 | 15 | 49 | 19 | 61 |
| 8 | 32 | 12 | 39 | 16 | 64 | 20 | 65 |

List 4.1

The Proof of Theorem 4.1:

Let

$$\begin{aligned} J_0 &= T_{s,1}(q) = \sum_p \int_0^1 \left| f(\alpha, p) \right|^{2(s-1)} \left| f(\alpha, q) \right|^2 d\alpha, \\ J_i &= \int_0^1 \left| f_i(\alpha) \right|^{2(s-1)} F_i(\alpha, q) d\alpha, \quad 1 \le i < k. \end{aligned} \tag{4.3}$$

Lemma 7 of paper [1] will be used in the following proofs, we restate here

**Lemma 4.1.**

$$\begin{aligned} J_0 &\le ZPS_{s-1}(P) + J_1, \\ J_i &\ll U_i + V_i, \\ U_i &= S_{s-1}(P_i)^{1/2} Z_{i+1}^{(2s-3)/2} \left( P(\tilde{H}_i \tilde{Z}_i)^2 Z_{i+1} S_{s-1}(P_{i+1}) \right)^{1/2} \\ V_i &= S_{s-1}(P_i)^{1/2} Z_{i+1}^{(2s-3)/2} (\tilde{H}_i \tilde{Z}_i J_{i+1})^{1/2}, \qquad 1 \le i < k, \end{aligned} \tag{4.4}$$

Where $\tilde{H}_i = \prod_{j \le i} H_j, \tilde{Z}_i = \prod_{j \le i} Z_j$ .

Let

$$U_{i-1} = V_{i-1} \cdot H_i^{-\alpha/2} P^{-\tau/2}, \quad |\alpha| < 1, \tau \geq 0, \alpha + \tau < 1, 1 \leq i \leq k. \tag{4.5}$$

The parameters $\alpha, \tau$ will be decided later.

Hence, it has

$$\left( P(\tilde{H}_{i-1}\tilde{Z}_{i-1})^2 Z_i S_{s-1}(P_i) \right) = (\tilde{H}_{i-1}\tilde{Z}_{i-1} J_i) \cdot H_i^{-\alpha} P^{-\tau} .$$

And

$$J_i = P(\tilde{H}_{i-1}\tilde{Z}_{i-1}) Z_i S_{s-1}(P_i) H_i^{\alpha} P^{\tau} . \tag{4.6}$$

On the other hand,

$$J_i \approx U_i \cdot H_{i+1}^{\alpha/2} P^{\tau/2} = S_{s-1}(P_i)^{1/2} Z_{i+1}^{(2s-3)/2} \left( P(\tilde{H}_i \tilde{Z}_i)^2 Z_{i+1} S_{s-1}(P_{i+1}) \right)^{1/2} H_{i+1}^{\alpha/2} P^{\tau/2} \tag{4.7}$$

Combine the two equalities, it has

$$P(\tilde{H}_{i-1}\tilde{Z}_{i-1}) Z_i S_{s-1}(P_i) H_i^{\alpha} P^{\tau} = S_{s-1}(P_i)^{1/2} Z_{i+1}^{(2s-3)/2} \left( P(\tilde{H}_i \tilde{Z}_i)^2 Z_{i+1} S_{s-1}(P_{i+1}) \right)^{1/2} H_{i+1}^{\alpha/2} P^{\tau/2}$$

It follows,

$$\frac{S_{s-1}(P_i)^{1/2}}{S_{s-1}(P_{i+1})^{1/2}} = \frac{Z_{i+1}^{(2s-2)/2} (H_i) H_{i+1}^{\alpha/2}}{P^{1/2} H_i^{\alpha} P^{\tau/2}} ,$$

and

$$\frac{S_{s-1}(P_i)}{S_{s-1}(P_{i+1})} = \frac{Z_{i+1}^{(2s-2)} (H_i)^{2-2\alpha} H_{i+1}^{\alpha}}{P^{1+\tau}} . \tag{4.8}$$

Let $S_{s-1}(X) = X^{\lambda_{s-1}}$, it has

$$\frac{P_i^{\lambda_{s-1}}}{P_{i+1}^{\lambda_{s-1}}} = Z_{i+1}^{(2s-2)} \frac{(P/Z_i^k)^{2-2\alpha} (P/Z_{i+1}^k)^{\alpha}}{P^{1+\tau}} .$$

i.e.

$$(Z_{i+1})^{\lambda_{s-1}} = Z_{i+1}^{(2s-2)} P^{1-\alpha-\tau} Z_i^{-k(2-2\alpha)} Z_{i+1}^{-k\alpha} .$$

and

$$\theta_i = \theta_{i+1} \frac{((2s-2) - \lambda_{s-1} - k\alpha)}{k(2-2\alpha)} + \frac{1}{2k} \frac{1-\alpha-\tau}{(1-\alpha)} . \tag{4.9}$$

Denote by $a = \dfrac{\left((2s-2) - \alpha k - \lambda_{s-1}\right)}{k(2-2\alpha)}, \ b = \dfrac{1}{2k} \dfrac{1-\alpha-\tau}{1-\alpha}$, it has

$$\theta_i = a^j \theta_{i+j} + b \frac{1-a^j}{1-a}, \text{ and } \theta_{k-i} = a^i \theta_k + b \frac{1-a^i}{1-a} . \tag{4.10}$$

Moreover, by (4.5) with $i = k$, it has

$$\left(P(\tilde{H}_{k-1}\tilde{Z}_{k-1})^2 Z_k S_{s-1}(P_k)\right)^{1/2} = (\tilde{H}_{k-1}\tilde{Z}_{k-1}J_k)^{1/2} \cdot H_k^{-\alpha/2}P^{-\tau/2}.$$

And it is easy to know that $J_k = P\tilde{H}_k\tilde{Z}_k S_{s-1}(P_k)$, hence

$$\left(P(\tilde{H}_{k-1}\tilde{Z}_{k-1})^2 Z_k S_{s-1}(P_k)\right) = (\tilde{H}_{k-1}\tilde{Z}_{k-1}P\tilde{H}_k\tilde{Z}_k S_{s-1}(P_k)) \cdot H_k^{-\alpha}P^{-\tau}$$

$$\Rightarrow 1 = H_k^{1-\alpha}P^{-\tau} \Rightarrow P^{\tau} = H_k^{1-\alpha} = (P/Z_k^k)^{1-\alpha}$$

$$\Rightarrow Z_k^{k(1-\alpha)} = P^{1-\alpha-\tau} \Rightarrow \theta_k = \frac{1}{k}\cdot\frac{1-\alpha-\tau}{(1-\alpha)}.$$

Let $\Delta(s)$ defined as before, there is

$$\theta_{k-i} = a^i\theta_k + b\frac{1-a^i}{1-a} = \frac{1-\alpha-\tau}{k(1-\alpha)+\Delta(s-1)} + \left(\frac{1-\alpha-\tau}{k(1-\alpha)} - \frac{1-\alpha-\tau}{k(1-\alpha)+\Delta(s-1)}\right)a^i.$$

Especially,

$$\theta = \theta_1 = \frac{1-\alpha-\tau}{k(1-\alpha)+\Delta(s-1)} + \left(\frac{1-\alpha-\tau}{k(1-\alpha)} - \frac{1-\alpha-\tau}{k(1-\alpha)+\Delta(s-1)}\right)a^{k-1}. \qquad (4.11)$$

On the other hand, by Lemma 2.1 with $i = 1$, it has

$$S_s(\tilde{P}) \ll Z^{2s-1}J_0(P) \approx Z^{2s-1}U_0H^{\alpha/2}P^{\tau/2} = Z^{2s-1}PZS_{s-1}(P)H^{\alpha/2}P^{\tau/2}.$$

And

$$\lambda_s = \frac{\lambda_{s-1}}{(1+\theta)} + \frac{(2s-k\alpha/2)\theta}{(1+\theta)} + \frac{(1+\alpha/2+\tau/2)}{(1+\theta)}.$$

Or,

$$\Delta(s) = \frac{\Delta(s-1)}{(1+\theta)} + \frac{(1-\alpha/2)(k\theta-1)}{(1+\theta)} + \frac{\tau}{2(1+\theta)}. \qquad (4.12)$$

When $k$ is greater, it may has

$$\theta \doteq \frac{(1-\alpha-\tau)}{(1-\alpha)k+\Delta(s-1)}. \qquad (4.13)$$

And let

$$(1-\alpha) = \frac{\beta\cdot\Delta(s-1)}{k}, \quad (1-\alpha-\tau) = \frac{\omega\cdot\Delta(s-1)}{k}. \qquad (4.14)$$

Substituting (4.13) and (4.14) in (4.12), it follows

$$\Delta(s) = \Delta(s-1)\left(1-\frac{3\omega}{2(k(\beta+1)+\omega)}\right) - \frac{\rho k}{2(k(\beta+1)+\omega)}. \qquad (4.15)$$

where $\rho = (1+\beta-\omega)$. And then

$$\Delta(s)=\Delta(s-i)\left(1-\frac{3\omega}{2(k(1+\beta)+\omega)}\right)^i-\frac{\rho k}{3\omega}\left(1-\left(1-\frac{3\omega}{2(k(1+\beta)+\omega)}\right)^i\right).$$

For $\Delta(2)=k-2$, hence

$$\Delta(i+2)=\left(k-2+\frac{\rho k}{3\omega}\right)\left(1-\frac{3\omega}{2(k(1+\beta)+\omega)}\right)^i-\frac{\rho k}{3\omega}, \tag{4.16}$$

Let

$$\Im(\beta,\omega)=\log\left(\frac{\rho k}{3\omega(k-2)+\rho k}\right)\Big/\log\left(1-\frac{3\omega}{2(k(1+\beta)+\omega)}\right). \tag{4.17}$$

From (4.16), we can know that $\Delta(i+2)$ will approach zero as $i$ tends to $\Im(\beta,\omega)$. Hence it has

**Lemma 4.2.** For sufficient large $k$, and arbitrary small $\varepsilon>0$, there is $s\le\Im(\beta,\omega)+2$ such that $\Delta(s)<\varepsilon.$

Clearly, for greater $k$, there is

$$\Im(\beta,\omega)\sim\frac{2k}{3}\left(1+\frac{\rho}{\omega}\right)\log\left(1+\frac{3\omega}{\rho}\right). \tag{4.18}$$

Denote by $x=\omega/\rho$, and let $\vartheta=0.60565351$, which is a root of the equation

$$\frac{3(1+x)}{1+3x}-\frac{1}{x}\log(1+3x)=0.$$

Hence,

$$\Im(\beta,\omega)\to\frac{2k}{3}(1+1/\vartheta)\log(1+3\vartheta)=1.83043k, \text{ as } x\to\vartheta.$$

We know that

$$G(k)\le 3+2u+2\left[\frac{\Delta(u)}{2\hat{\sigma}}\right].$$

where $\hat{\sigma}=\frac{\log(1+1/k)}{4(1+\lambda)}$, $\lambda=\log k+\log\log k+O\left(\frac{\log\log k}{\log k}\right)$.

i.e.

$$G(k)\le 3+2u+(4+o(1))\Delta(u)k\log k.$$

Take $u=\Im(\beta,\omega)+2,$ and let $x\to\vartheta$, Theorem 4.1 is proved.

The Proof of Theorem 4.2:

For smaller $k$, $\Delta(s)$ may be followed by recursions (4.11) and (4.12) and the initial $\Delta(2)$ by choosing optimally values $\alpha$ and $\tau$ in turn. With the restriction of the ability of PC, we have lessened the search range of parameters $\alpha$, $\tau$ only in five digits. The results in List 4.2 are obtained by PC. A completed list including intermediate results is posted behind as appendix 1.

| $k$ | $s$ | $\Delta(s)$ | $k$ | $s$ | $\Delta(s)$ |
|---|---|---|---|---|---|
| 5 | 8 | 0.000000 | 13 | 21 | 0.000000 |
| 6 | 10 | 0.000000 | 14 | 22 | 0.000000 |
| 7 | 11 | 0.000000 | 15 | 24 | 0.000000 |
| 8 | 13 | 0.000000 | 16 | 25 | 0.000000 |
| 9 | 14 | 0.000000 | 17 | 27 | 0.000000 |
| 10 | 16 | 0.000000 | 18 | 28 | 0.000000 |
| 11 | 18 | 0.000000 | 19 | 30 | 0.000000 |
| 12 | 19 | 0.000000 | 20 | 32 | 0.000000 |

List 4.2

As the proof of Theorem 1.2, with (3.2) and List 4.2, we take $v(k)$ as in the following list

| $k$ | 5 | 6 | 7 | 8 | 9 | 10 | 11 | 12 | 13 | 14 | 15 | 16 | 17 | 18 | 19 | 20 |
|---|---|---|---|---|---|---|---|---|---|---|---|---|---|---|---|---|
| $v(k)$ | 1 | 1 | 1 | 6 | 1 | 1 | 1 | 1 | 1 | 1 | 1 | 14 | 1 | 1 | 1 | 1 |

List 4.3

And Theorem 4.2 is followed.

**5. Further Improvements (2)**

In this section, we will present further improvement when $k$ is larger.

**Theorem 5.1.** For sufficient large $k$,

$$G(k) \le \begin{cases} 3.182k + o(k), & \text{if } k \text{ is not a power of } 2, \\ 4k, & \text{otherwise.} \end{cases} \tag{5.1}$$

The Proof of Theorems 5.1:

The notations and symbols used here will be same as before. Define

$$F_i^{(2)}(\alpha,q)=\sum_{1\le h_1\le H_1}\cdots\sum_{1\le h_i\le H_i}\sum_{p_1\in\mathscr{P}_1}\cdots\sum_{p_i\in\mathscr{P}_i}\sum_{x,y}e((\Psi_i(x)-\Psi_i(y))\alpha q^k), \tag{5.2}$$

$$\widehat{F}_i^{(2)}(\alpha,q)=\sum_{1\le h_1\le H_1}\cdots\sum_{1\le h_i\le H_i}\sum_{p_1\in\mathscr{P}_1}\cdots\sum_{p_i\in\mathscr{P}_i}\sum_{x\ne y}e((\Psi_i(x)-\Psi_i(y))\alpha q^k),\quad 1\le i\le k. \tag{5.3}$$

Let

$$\begin{aligned} J_0&=\sum_p\int_0^1|f(\alpha,p)|^{2(s-1)}|f(\alpha,q)|^2\,d\alpha,\\ J_i&=\int_0^1|f_i(\alpha)|^{2(s-1)}F_i(\alpha,q)d\alpha,\qquad 1\le i<k. \end{aligned} \tag{5.4}$$

There is

**Lemma 5.1.**

$$\begin{aligned} &J_0\le ZPS_{s-1}(P)+J_1,\\ &J_i\ll U_i+V_i,\\ &U_i=P^{1/2}\tilde{H}_i\tilde{Z}_iS_{s-1}(P_i),\quad V_i=S_{s-1}(P_i)^{1/2}Z_{i+1}^{(2s-3)/2}(\tilde{H}_i\tilde{Z}_iJ_{i+1})^{1/2},\qquad 1\le i<k, \end{aligned} \tag{5.5}$$

where $\tilde{H}_i=\prod_{j\le i}H_j$, $\tilde{Z}_i=\prod_{j\le i}Z_j$ .

Proof. By Cauchy inequality,

$$J_i=\int_0^1|f_i(\alpha)|^{2(s-1)}F_i(\alpha,q)d\alpha\le\left(\int_0^1|f_i(\alpha)|^{2(s-1)}d\alpha\right)^{1/2}\left(\int_0^1|f_i(\alpha)|^{2(s-1)}|F_i(\alpha,q)|^2\,d\alpha\right)^{1/2}$$

$$\le\left(\int_0^1|f_i(\alpha)|^{2(s-1)}d\alpha\right)^{1/2}\left(\tilde{H}_i\tilde{Z}_i\int_0^1|f_i(\alpha)|^{2(s-1)}F_i^{(2)}(\alpha,q)d\alpha\right)^{1/2}$$

$$\le\left(\int_0^1|f_i(\alpha)|^{2(s-1)}d\alpha\right)^{1/2}(\tilde{H}_i\tilde{Z}_i)^{1/2}\left(\int_0^1|f_i(\alpha)|^{2(s-1)}P(\tilde{H}_i\tilde{Z}_i)d\alpha+\int_0^1|f_i(\alpha)|^{2(s-1)}\widehat{F}_i(\alpha,q)^{(2)}d\alpha\right)^{1/2}$$

$$\le S_{s-1}(P_i)^{1/2}(\tilde{H}_i\tilde{Z}_i)^{1/2}\left(P(\tilde{H}_i\tilde{Z}_i)S_{s-1}(P_i)+\int_0^1\left|\sum_{p\in\mathscr{P}_{i+1}}f_{i+1}(\alpha,p)\right|^{2(s-1)}\widehat{F}_i(\alpha,q)^{(2)}d\alpha\right)^{1/2}$$

$$\le S_{s-1}(P_i)^{1/2}(\tilde{H}_i\tilde{Z}_i)^{1/2}\left(P(\tilde{H}_i\tilde{Z}_i)S_{s-1}(P_i)+Z_{i+1}^{2s-3}\int_0^1\sum_{p\in\mathscr{P}_{i+1}}|f_{i+1}(\alpha,p)|^{2(s-1)}\widehat{F}_i(\alpha,q)^{(2)}d\alpha\right)^{1/2}$$

$$\ll P^{1/2}\tilde{H}_i\tilde{Z}_iS_{s-1}(P_i)+S_{s-1}(P_i)^{1/2}(Z_{i+1}^{2s-3})^{1/2}(\tilde{H}_i\tilde{Z}_i)^{1/2}\left(\int_0^1|f_{i+1}(\alpha)|^{2(s-1)}F_{i+1}(\alpha,q)d\alpha\right)^{1/2}$$

$$\ll P^{1/2}\tilde{H}_i\tilde{Z}_iS_{s-1}(P_i)+S_{s-1}(P_i)^{1/2}(Z_{i+1}^{2s-3})^{1/2}(\tilde{H}_i\tilde{Z}_i)^{1/2}(J_{i+1})^{1/2}$$

□

Let

$$U_{i-1} = V_{i-1} \cdot (H_i)^{-\alpha/2} P^{-\tau/2}, \quad 1 \le i < k. \tag{5.6}$$

The parameters $\alpha, \tau$ will be determined later. So, it follows

$$J_i = \frac{P(\tilde{H}_{i-1}\tilde{Z}_{i-1})S_{s-1}(P_{i-1})}{Z_i^{(2s-3)}}(H_i)^{\alpha} P^{\tau}. \tag{5.7}$$

Besides,

$$J_i \approx U_i (H_{i+1})^{\alpha/2} P^{\tau/2} = P^{1/2}\tilde{H}_i\tilde{Z}_i S_{s-1}(P_i)(H_{i+1})^{\alpha/2} P^{\tau/2}. \tag{5.8}$$

Combine the two identities above, it follows

$$\frac{S_{s-1}(P_{i-1})}{S_{s-1}(P_i)} = Z_i^{(2s-2)} \frac{P^{-1/2} H_i (H_{i+1})^{\alpha/2}}{(H_i)^{\alpha} P^{\tau/2}}. \tag{5.9}$$

i.e.

$$(Z_i)^{\lambda_{s-1}} = Z_i^{(2s-2)} P^{-(\tau+1)/2} (P / Z_i^k)^{1-\alpha} (P / Z_{i+1}^k)^{\alpha/2}$$

Or,

$$P^{\theta_i(\lambda_{s-1}-(2s-2)+k(1-\alpha))} = P^{(1-\alpha-\tau)/2} P^{-\theta_{i+1}k\alpha/2}.$$

And,

$$\theta_{i+1} = \theta_i \frac{2((2s-2)-\lambda_{s-1}-k(1-\alpha))}{k\alpha} + \frac{(1-\alpha-\tau)}{k\alpha}. \tag{5.10}$$

Denote by $a = 2\dfrac{\alpha k - \Delta(s-1)}{\alpha k}$, $b = \dfrac{(1-\alpha-\tau)}{\alpha k}$.

Then it has

$$\theta_{i+j} = a^j \theta_i + b\frac{1-a^j}{1-a}, \text{ and } \theta_k = a^i \theta_{k-i} + b\frac{1-a^i}{1-a}.$$

Moreover, by (5.6) with $i = k$, it has

$$P(\tilde{H}_{k-1}\tilde{Z}_{k-1})^2 S_{s-1}(P_{k-1})^2 = S_{s-1}(P_{k-1}) Z_k^{(2s-3)} (\tilde{H}_{k-1}\tilde{Z}_{k-1}J_k)(H_k)^{-\alpha} P^{-\tau}.$$

i.e.

$$P(\tilde{H}_{k-1}\tilde{Z}_{k-1})S_{s-1}(P_{k-1}) = Z_k^{(2s-3)} J_k (H_k)^{-\alpha} P^{-\tau}.$$

It is easy to know that $J_k = P\tilde{H}_k\tilde{Z}_k S_{s-1}(P_k)$, hence

$$P(\tilde{H}_{k-1}\tilde{Z}_{k-1})S_{s-1}(P_{k-1}) = Z_k^{(2s-3)} P\tilde{H}_k\tilde{Z}_k S_{s-1}(P_k)(H_k)^{-\alpha} P^{-\tau}.$$

And

$$\frac{S_{s-1}(P_{k-1})}{S_{s-1}(P_k)} = Z_k^{(2s-2)} H_k (H_k)^{-\alpha} P^{-\tau}.$$

i.e.

$$(Z_k)^{\lambda_{s-1}-(2s-2)+k(1-\alpha)} = P^{1-\alpha-\tau}$$

It follows

$$\theta_k = \frac{1-\alpha-\tau}{\Delta(s-1)-k\alpha}. \tag{5.11}$$

Hence,

$$\theta_{k-i} = \frac{(1-\alpha-\tau)}{2\Delta(s-1)-\alpha k} + \frac{(1-\alpha-\tau)\Delta(s-1)}{(\Delta(s-1)-k\alpha)(2\Delta(s-1)-\alpha k)} a^{-i}.$$

Especially,

$$\theta = \theta_1 = \frac{(1-\alpha-\tau)}{2\Delta(s-1)-\alpha k} + \frac{(1-\alpha-\tau)\Delta(s-1)}{(\Delta(s-1)-k\alpha)(2\Delta(s-1)-\alpha k)} a^{-k+1}. \tag{5.12}$$

Furthermore, by Lemma 2.1 with $i=1$, it has

$$S_s(\tilde{P}) \ll Z^{2s-1}J_0(P) \approx Z^{2s-1}U_0 H^{\alpha/2}P^{\tau/2} = Z^{2s-1}PZS_{s-1}(P)H^{\alpha/2}P^{\tau/2}.$$

And

$$\lambda_s = \frac{\lambda_{s-1}}{(1+\theta)} + \frac{(2s-k\alpha/2)\theta}{(1+\theta)} + \frac{(1+\alpha/2+\tau/2)}{(1+\theta)}.$$

i.e.

$$\Delta(s) = \frac{\Delta(s-1)}{(1+\theta)} + \frac{(1-\alpha/2)(k\theta-1)}{(1+\theta)} + \frac{\tau}{2(1+\theta)}. \tag{5.13}$$

Let

$$\alpha = \frac{\beta\cdot\Delta(s-1)}{k}, \quad (1-\alpha-\tau) = \frac{\omega\cdot\Delta(s-1)}{k}. \tag{5.14}$$

Then it follows

$$\Delta(s) = \Delta(s-1)\left(1-\frac{4\omega+\delta(2+\beta)}{2(k(2-\beta)+\delta+\omega)}\right) + \frac{k(\omega+2\delta-\rho)}{2(k(2-\beta)+\delta+\omega)}, \tag{5.15}$$

where $\rho = 2-(\omega+\beta)$, $\delta = \frac{\omega}{(1-\beta)}\left(2\frac{\beta-1}{\beta}\right)^{-k+1}$.

And then

$$\Delta(s) = \Delta(s-i)\left(1-\frac{4\omega+\delta(2+\beta)}{2(k(2-\beta)+\delta+\omega)}\right)^i + \frac{k(\omega+2\delta-\rho)}{4\omega+\delta(2+\beta)}\left(1-\left(1-\frac{4\omega+\delta(2+\beta)}{2(k(2-\beta)+\delta+\omega)}\right)^i\right).$$

i.e.

$$\Delta(d+i) = \left(\Delta(d)+\frac{k(\rho-\omega-2\delta)}{4\omega+\delta(2+\beta)}\right)\left(1-\frac{4\omega+\delta(2+\beta)}{2(k(2-\beta)+\delta+\omega)}\right)^i - \frac{k(\rho-\omega-2\delta)}{4\omega+\delta(2+\beta)}. \tag{5.16}$$

Let

$$\Im(\beta,\omega)=\log\left(\frac{k(\rho-\omega-2\delta)}{(4\omega+\delta(2+\beta))\Delta(d)+k(\rho-\omega-2\delta)}\right)\Big/\log\left(1-\frac{4\omega+\delta(2+\beta)}{2(k(2-\beta)+\delta+\omega)}\right). \tag{5.17}$$

From (5.16), we can know that $\Delta(i+d)$ will approach zero as $i$ tends to $\Im(\beta,\omega)$. Hence it has

**Lemma 5.2.** For sufficient large $k$, and arbitrary small $\varepsilon>0$, there is $s\le\Im(\beta,\omega)+d$ such that $\Delta(s)<\varepsilon.$

It is clear that, when $k$ is greater,

$$\Im(\beta,\omega)\sim\frac{2k(2-\beta)}{4\omega+\delta(2+\beta)}\log\left(1+\frac{4\omega+\delta(2+\beta)}{(\rho-\omega-2\delta)}\frac{\Delta(d)}{k}\right). \tag{5.18}$$

We choose $\Delta(d)$ such that $\Delta(d)\le k/3,$ then we can take $\beta=-3,$ and $0<\omega<\varepsilon$. Hence, it has $\delta=o(\omega).$ And then let $\omega\to 0$, by (5.18), it follows

$$\Im(\beta,\omega)\sim 2\Delta(d). \tag{5.19}$$

On the other hand, from § 4 (there take $\beta=1,\ \tau=0$), we can know that

$$d\le\frac{4k}{3}\log 2=0.9242k\,. \tag{5.20}$$

We have known that (see [3])

$$G(k)\le 3+2u+(4+o(1))\Delta(u)k\log k.$$

Take $u=d+\Im(\beta,\omega),$ Theorem 5.1 is proved.

In fact, Theorem 5.1 can be also proved along the way of section 4 as following.

The Second Proof of Theorem 5.1:

By (4.11) and (4.12), there is

$$\Delta(s)=\Delta(s-1)\left(1-\frac{3\omega+(2-\beta)\delta}{2(k(\beta+1)+\delta+\omega)}\right)-\frac{k(\rho-\delta)}{2(k(\beta+1)+\delta+\omega)}.$$

Where $1-\alpha=\dfrac{\beta\cdot\Delta(s-1)}{k},1-\alpha-\tau=\dfrac{\omega\cdot\Delta(s-1)}{k},\delta=\dfrac{\omega}{\beta}\left(\dfrac{\beta-1}{2\beta}\right)^{k-1},\ \rho=1+\beta-\omega.$

And then,

$$\Delta(d+i)=\left(\Delta(d)+\frac{k(\rho-\delta)}{3\omega+(2-\beta)\delta}\right)\left(1-\frac{3\omega+(2-\beta)\delta}{2(k(\beta+1)+\delta+\omega)}\right)^i-\frac{k(\rho-\delta)}{3\omega+(2-\beta)\delta}. \qquad (5.21)$$

Let

$$\Im(\beta,\omega)=\log\left(\frac{k(\rho-\delta)}{\Delta(d)(3\omega+(2-\beta)\delta)+k(\rho-\delta)}\right)\Big/\log\left(1-\frac{3\omega+(2-\beta)\delta}{2(k(\beta+1)+\delta+\omega)}\right). \qquad (5.22)$$

From (5.21), we can know that $\Delta(i+d)$ will approach zero as $i$ tends to $\Im(\beta,\omega)$.

Clearly, when $k$ is greater,

$$\Im(\beta,\omega)\sim\frac{2k(\beta+1)}{3\omega+(2-\beta)\delta}\log\left(1+\frac{3\omega+(2-\beta)\delta}{(\rho-\delta)}\frac{\Delta(d)}{k}\right). \qquad (5.23)$$

We choose $\Delta(d)$ such that $\Delta(d)\le k/3$, then it can be taken that $\beta=6$, and $\delta=o(\omega)$. Take

$0<\omega<\varepsilon$, and then let $\varepsilon\to 0$, by (5.23), it follows

$$\Im(\beta,\omega)\sim 2\Delta(d).$$

The rest is same as the first proof.

**6. The treat for the singular cases**

In this section, we supplement an argument for the singular cases which is saved in the proofs of Lemmas 4.1 and 5.1.

Let $\Phi(x)=\Phi(x;h_1,\cdots,h_r;p_1,\cdots,p_r)=\Delta_r(x^k;h_1,\cdots,h_r;p_1,\cdots,p_r)$, and denote by $\mathcal{N}_r(\Phi)$ the number of solutions of the equation

$$x_1^k+\cdots+x_{s-1}^k+\Phi(x)=y_1^k+\cdots+y_{s-1}^k+\Phi(y),\quad x_i,y_i\in\mathscr{C}_r, 1\le i<s,\ x,y\in P. \qquad (6.1)$$

Let $\mathcal{N}_r=\sum_{h_1,\cdots,h_r}\sum_{p_1,\cdots,p_r}\mathcal{N}_r(\Phi)$.

From (4.1), we can know that, provided that $U_i\ll V_i$,

$$\mathcal{N}_r\ll\Lambda\cdot\Gamma,$$
$$\Lambda=\prod_{r<t\le k}Z_t^{(2s-3)/2^{t-r-1}}\prod_{r<t<k}(\tilde{Z}_tS_{s-1}(P_t))^{1/2^{t-r}}(\tilde{Z}_kS_{s-1}(P_k))^{1/2^{k-r-1}}, \qquad (6.2)$$
$$\Gamma=\prod_{r<t<k}(\tilde{H}_t)^{1/2^{t-r}}(\tilde{H}_k\cdot P)^{1/2^{k-r-1}}.$$

Note that $\Lambda$ is only a function of $Z_i$ 's and $S_s(P_i)$ 's.

For an integer coefficient polynomial $f(x)$, denote by $c_f$ the greatest common divisor of the coefficients of $f(x)$.

For a number set $A$, $\varpi(A)$ denote the number of distinct elements in set $A \mod p^k$.

Suppose that $p^a \| c_{\Phi'}, a \geq 0$. If $a = 0$, i.e. $(c_{\Phi'}, p) = 1$, let $\{x_l\}_1^d$ be the set of solutions of equation $\Phi'(x) \equiv 0 \mod p$, and for a $x_l$, denote by

$$X^{(l)} = \{x \in P \mid x \equiv x_l \mod p\}, \quad X_n = \{x \in X^{(l)} \mid \Phi(x) \equiv n \mod p^k\},$$
$$\mathfrak{A}_i = \bigcup_{p^i \leq \varpi(X_n) < p^{i+1}} X_n.$$

We divide $\mathfrak{A}_i$ into subsets $\{\mathfrak{A}_{i,j}\}_1^{\mu}$ such that in each subset $\mathfrak{A}_{i,j}$ equation $\Phi(x) \equiv n \mod p^k$ has at most a solution. Denote by $Q = \max_j \{ |\mathfrak{A}_{i,j}| \}$, then it is easy to know that $p^i \leq \mu < p^{i+1}$ and $Q \leq |X^{(l)}| / p^i \leq P / p^{i+1}$. By Cauchy inequality, it has

$$\left| \sum_{x \in X^{(l)}} e(\Phi(x)\alpha) \right|^2 = \left| \sum_i \sum_{x \in \mathfrak{A}_i} e(\Phi(x)\alpha) \right|^2 \leq k \sum_i \left| \sum_{x \in \mathfrak{A}_i} e(\Phi(x)\alpha) \right|^2 \leq k\mu \sum_i \sum_{1 \leq j \leq \mu} \left| \sum_{x \in \mathfrak{A}_{i,j}} e(\Phi(x)\alpha) \right|^2$$

Now we can apply the iterative differencing method to the each subset $\mathfrak{A}_{i,j}$ as in the normal case, but $H_t|_{singular} \leq Q / Z_t^k \leq P / (p^{i+1} Z_t^k) = H_t|_{normal} \cdot p^{-(i+1)}$.

For arbitrary non-negative integer $v$, it is easy to know that $I_v = \sum_{0 \leq i \leq v} i \cdot 2^{-i} + (v+2) \cdot 2^{-v} = 2$.

Denote by $\tilde{H}_{r,t} = \prod_{r < i \leq t} H_i$, $\tilde{Z}_{r,t} = \prod_{r < i \leq t} Z_i$, then,

$$\begin{aligned} &\mu \cdot \mu \cdot \prod_{r<t<k} (\tilde{H}_{r,t}|_{singular})^{1/2^{t-r}} (\tilde{H}_{r,k}|_{singular} \cdot Q)^{1/2^{k-r-1}} \\ &\leq \prod_{r<t<k} (\tilde{H}_{r,t}|_{normal})^{1/2^{t-r}} (\tilde{H}_{r,k}|_{normal} \cdot P)^{1/2^{k-r-1}} \cdot p^{(i+1)(2-I_{k-r-1})} \quad . \qquad (6.3) \\ &\leq \prod_{r<t<k} (\tilde{H}_{r,t}|_{normal})^{1/2^{t-r}} (\tilde{H}_{r,k}|_{normal} \cdot P)^{1/2^{k-r-1}} \end{aligned}$$

With (6.2), this indicates that the contribution from the singular parts no more than the one from the normal parts times a constant.

For the case $0 < a < k$, consider the equation $\Phi(x) \equiv n \mod p^k$, as it can be changed as $\Phi(x) - \Phi(0) \equiv n - \Phi(0) \mod p^k$, so for simplicity, assume the constant term $\Phi(0) = 0$. Hence $p^a \mid \Phi(x)$, and the solvability of equation $\Phi(x) \equiv n \mod p^k$ implies $p^a \mid n$. Upon that the equation will be degraded as $\bar{\Phi}(x) \equiv \bar{n} \mod p^{k-a}$, where $\bar{\Phi}(x) = \Phi(x) / p^a$, $\bar{n} = n / p^a$.

Let $\Theta = \left\{ x \in P \mid \bar{\Phi}'(x) \equiv 0 \mod p \right\}$, and $\Omega = P \setminus \Theta$.

We know that in set $\Omega$ the equation $\bar{\Phi}(x) \equiv \bar{n} \mod p^{k-a}$ has at most finite solutions. As in the normal case, it will bring about $H_{r+1}\big|_{\Omega}$ differencing functions of $\bar{\Phi}(x)$ but with step $p^{k-a}$, where

$$H_{r+1}\big|_{\Omega} \le P / p^{k-a} \le H_{r+1}\big|_{normal} \cdot p^{a}. \tag{6.4}$$

On the other hand, write $\Phi(x) = \Delta_r(x^k, h_1, \cdots, h_r)$, by the expanding formula of difference (see [3]), there is $p^a \mid h_1 \cdots h_r$, if $p^a \mid c_{\Phi}$, $p > k$. This means

$$\prod_{1 \le i \le r} H_i\big|_{singular} \ll \left( \prod_{1 \le i \le r} H_i\big|_{normal} \right) \cdot p^{-a}. \tag{6.5}$$

Roughly speaking, the singular cases are distributed in proportion. Hence,

$$\prod_{1 \le i \le r+1} H_i\big|_{singular} \ll \prod_{1 \le i \le r+1} H_i\big|_{normal}. \tag{6.6}$$

For the set $\Theta$, let $\left\{ x_l \right\}_1^d$ be the set of solutions of equation $\bar{\Phi}'(x) \equiv 0 \mod p$, and for a $x_l$, denote by

$$X^{(l)} = \left\{ x \in P \mid x \equiv x_l \mod p \right\}, \quad \bar{X}_n = \left\{ x \in X^{(l)} \mid \bar{\Phi}(x) \equiv \bar{n} \mod p^{k-a} \right\},$$

$$\bar{\mathfrak{A}}_i = \bigcup_{p^i \le \bar{\varpi}(\bar{X}_n) < p^{i+1}} \bar{X}_n.$$

Where $\bar{\varpi}(A)$ denote the number of distinct elements in the number set $A \mod p^{k-a}$.

Divide $\bar{\mathfrak{A}}_i$ into subsets $\left\{ \bar{\mathfrak{A}}_{i,j} \right\}_1^{\mu}$ such that in each subset $\bar{\mathfrak{A}}_{i,j}$ equation $\bar{\Phi}(x) \equiv \bar{n} \mod p^{k-a}$ has at most a solution. Then

$$H_{r+1}\big|_{singular} \le Q / p^{k-a} \le P / (p^{i+1} p^{k-a}) = H_{r+1}\big|_{normal} \cdot p^{a-(i+1)}.$$

Thus,

$$\begin{aligned} &\mu \cdot \mu \cdot \prod_{r<t<k} (\tilde{H}_{r,t}\,|_{singular})^{1/2^{t-r}} (\tilde{H}_{r,k}\,|_{singular} \cdot Q)^{1/2^{k-r-1}} \\ &\le \prod_{r<t<k} (\tilde{H}_{r,t}\,|_{normal})^{1/2^{t-r}} (\tilde{H}_{r,k}\,|_{normal} \cdot P)^{1/2^{k-r-1}} \cdot p^{(i+1)(2-I_{k-r-1})+a} \\ &\le \prod_{r<t<k} (\tilde{H}_{r,t}\,|_{normal})^{1/2^{t-r}} (\tilde{H}_{r,k}\,|_{normal} \cdot P)^{1/2^{k-r-1}} \cdot p^{a} \end{aligned} \tag{6.7}$$

Together with (6.5), our assertion is followed.

For the case $a \ge k$, clearly, $H_{r+1}\big|_{singular} \ll P$, as all the possible differencing functions are of $\Phi(x+h) - \Phi(x)$, $|h| < P$. And by (6.5), the assertion holds.

Finally, it should be mentioned that estimation (6.2), which is deduced under a condition that $U_i \ll V_i$, is a main tool in the argument above and applied in the normal cases and the singular cases implicitly at the same time. It is not difficult to know that the condition in the singular cases is somewhat stronger than in the normal cases for $Q < P$ in general.

Without the condition, estimation (6.2) will be changed as

$$\begin{aligned} &\mathcal{N}_r \ll \sum_{r\le v<k} \Lambda_v \cdot \Gamma_v + \Lambda \cdot \Gamma, \\ &\Lambda_r = \tilde{Z}_r Z_{r+1}^{(2s-2)} S_{s-1}(P_{r+1}), \quad \Gamma_r = \tilde{H}_r P. \\ &\Lambda_v = \prod_{r<t\le v+1} Z_t^{(2s-3)/2^{t-r-1}} \prod_{r<t\le v} (\tilde{Z}_t S_{s-1}(P_t))^{1/2^{t-r}} \left(\tilde{Z}_{v+1} S_{s-1}(P_{v+1})\right)^{1/2^{v-r}}, \quad r<v<k. \\ &\Gamma_v = \prod_{r<t\le v} (\tilde{H}_t)^{1/2^{t-r}} \left(P\tilde{H}_v\right)^{1/2^{v-r}}, \quad r<v<k. \end{aligned} \tag{6.2'}$$

In the following, we will show that the postulate may be released, that is, we may verify our assertion by (6.2') instead of (6.2).

At first we note that if $x_l$ is only a simple root of $\Phi'(x) \mod p$, then there are actually

$$\left|\mathfrak{A}_i\right| \ll P/p^i, \text{ and } Q \ll P/p^{2i}. \tag{6.8}$$

So, our assertion may be followed by (6.2') as before.

If $x_l$ is a multiple root, the investigation will be a little more complex. Let

$$\Phi(x) - n \equiv A(x)\cdot B(x) \mod p^k, \tag{6.9}$$

where $A(x) = \prod_{0\le i\le m} (x - \vartheta_i)$, $\vartheta_i \equiv x_l \mod p$, $0 \le i \le m$. $B(x_l) \not\equiv 0 \mod p$.

Suppose that $p^{\alpha_i} \,\|\, \left(\vartheta_0 - \vartheta_i\right), 1 \le i \le m,$ without loss the generality, assume $\alpha_1 \ge \alpha_2 \ge \cdots \ge \alpha_m$.

Then there are $\left\{\bar{\vartheta}_i\right\}_1^m$, which are the roots of $\Phi^{(i)}(x) \equiv 0 \mod p^k$, and $\vartheta_i \equiv \bar{\vartheta}_i \mod p^{\alpha_i}$, $1 \le i \le m$. Suppose that $X_n \subseteq \mathfrak{A}_t$, i.e. $p^t \le \varpi(X_n) < p^{t+1}$, then it is clear that

$$t \le \sum_{1\le i\le m} \alpha_i\,. \tag{6.10}$$

Besides, for any $x \in X_n$, it should be that $x \equiv \bar{\vartheta}_1 \mod p^{\alpha_1}$, and it is easy to know that actually $\mu = O(p^t)$, so it has

$$\left|\mathfrak{A}_t\right| \ll P/p^{\alpha_1}, \text{ and } Q \ll P/p^{\alpha_1+t}. \tag{6.11}$$

Moreover, with $m$ equations $\Phi^{(i)}(\overline{\vartheta}_i) \equiv 0 \mod p^{\alpha_i}, 1 \le i \le m$, which will be formed a algebraic system of the variables $h_1, \cdots, h_r$ of $\Phi(x) (= \Delta_r(x^k; h_1, \cdots, h_r))$, $(h_i, p) = 1, 1 \le i \le r$, $(p > k)$, and the sizes of variables $h_1, \cdots, h_r$ will be reduced to about $p^{-(\alpha_1 + \cdots + \alpha_m)}$ in total for a given array $(\overline{\vartheta}_1, \cdots, \overline{\vartheta}_m)$ ( by elimination method ), consequently

$$\prod_{1 \le i \le r} H_i \big|_{singular} \ll \left( \prod_{1 \le i \le r} H_i \big|_{normal} \right) \cdot p^{-(\alpha_1 + \cdots + \alpha_m)}. \qquad (6.12)$$

With (6.2'), (6.10) ~ (6.12), our assertion may be deduced as before.

It is easy to calculate that the Jacobian criterion for the functions independence of functions $\Phi^{(i)}(x; h_1, \cdots, h_r), 1 \le i \le r$, is equal to

$$\frac{\partial(\Phi^{(1)}, \cdots, \Phi^{(r)})}{\partial(h_1, \cdots, h_r)} = (-1)^r \left( \prod_{1 < i \le r+1} [k]_i \right) (h_1 \cdots h_r)^{r-1} \cdot M \cdot R(x, h_1, \cdots, h_r). \qquad (6.13)$$

Where $[k]_i = k(k-1)\cdots(k-i+1)$, $M = \prod_{i<j} (h_i - h_j)$, i.e. the van der Monde's determinant of variables $h_1, \cdots, h_r$, and $R(x, h_1, \cdots, h_r)$ is a homogenous polynomial of variables $h_1, \cdots, h_r, x$ of degree $r(k-2r)$, in which the highest term of variable $x$ is $c \cdot x^{(k-2r)r}$, the constant $c$ can be defined by the approach of taking the derivation of order $(k-2r)r$ in the both sides of (6.13) with respect to variable $x$, e.g. $c = C_{k-2}^2,\ 2C_{k-2}^2 C_{k-3}^2 C_{k-4}^2 / 3$ for $r = 2, 3$ respectively.

This is in favor of the elimination method can be processed properly.

There is assumed that $k \ge 2r$, the case $k < 2r$ is similar provided that the order of Jacobian criterion is changed from $r$ to $k-r$ for arbitrary $k-r$ ones of variables $h_1, \cdots, h_r$.

In addition, it may be worth noting that the singularity of $\Phi(x)$ will not only influence the size of $\tilde{H}_r$, but also $\tilde{H}_{r-1}, \tilde{H}_{r-2}, \cdots$, that is, there are so called $\tilde{H}_{r-i}\big|_{singular}, i = 0, 1, \cdots$. and it is easy to know that

$\tilde{H}_{r-i}\big|_{singular} \ll \tilde{H}_{r-i}\big|_{normal} \cdot p^{-\omega_i}, \omega_i = \sum_{i<j\le m} \alpha_j, i = 0, 1, \cdots, m-1.$ Hence, there is

$$\prod_{1 \le i \le r} (\tilde{H}_i\big|_{singular})^{1/2^i} (\tilde{H}_r\big|_{singular})^{1/2^r} \ll \prod_{1 \le i \le r} (\tilde{H}_i\big|_{normal})^{1/2^i} (\tilde{H}_r\big|_{normal})^{1/2^r} \cdot \prod_{1 \le i \le m} p^{-\alpha_i / 2^{r-i}}. \qquad (6.14)$$

Clearly, (6.14) is more accurate than (6.12), and indicates that the multiplicity $m$ has a stronger restriction for the number of polynomials $\Phi(x; h_1, \cdots, h_r)$ .

## Appendix Intermediate Results of Recursions for Theorem 4.2

| $k$ | $s$ | $\alpha$ | $\tau$ | $\theta$ | $\Delta(s)$ |
|---|---|---|---|---|---|
| 5 | 3 | 0 | 0 | 0.125120 | 2.333618 |
| | 4 | 0 | 0 | 0.136680 | 1.774482 |
| | 5 | -0.999900 | 1.999800 | 0.000009 | 1.274485 |
| | 6 | -0.999900 | 1.999800 | 0.000009 | 0.774495 |
| | 7 | -0.999900 | 1.999800 | 0.000009 | 0.274512 |
| | 8 | -0.995600 | 1.021900 | 0.095118 | 0.000000 |
| 6 | 3 | 0 | 0 | 0.100009 | 3.272749 |
| | 4 | 0 | 0 | 0.107879 | 2.635686 |
| | 5 | 0 | 0 | 0.115887 | 2.088927 |
| | 6 | -0.999900 | 1.999800 | 0.000007 | 1.588930 |
| | 7 | -0.999900 | 1.999800 | 0.000007 | 1.088938 |
| | 8 | -0.999900 | 1.999800 | 0.000008 | 0.588953 |
| | 9 | -0.999900 | 1.999800 | 0.000008 | 0.088973 |
| | 10 | -0.153400 | 0.179400 | 0.139009 | 0.000000 |
| 7 | 3 | 0 | 0 | 0.083334 | 4.230771 |
| | 4 | 0 | 0 | 0.089044 | 3.538957 |
| | 5 | 0 | 0 | 0.094897 | 2.925605 |
| | 6 | 0 | 0 | 0.100775 | 2.390162 |
| | 7 | -0.999900 | 1.999800 | 0.000006 | 1.890164 |
| | 8 | -0.999900 | 1.999800 | 0.000006 | 1.390172 |
| | 9 | -0.999900 | 1.999800 | 0.000007 | 0.890184 |
| | 10 | -0.999900 | 1.999800 | 0.000007 | 0.390202 |
| | 11 | -0.998600 | 1.521400 | 0.033196 | 0.000000 |
| 8 | 3 | 0 | 0 | 0.071429 | 5.200000 |
| | 4 | 0 | 0 | 0.075758 | 4.467607 |
| | 5 | 0 | 0 | 0.080209 | 3.804151 |
| | 6 | 0 | 0 | 0.084719 | 3.209961 |
| | 7 | 0 | 0 | 0.089214 | 2.684204 |
| | 8 | -0.999900 | 1.999800 | 0.000005 | 2.184207 |
| | 9 | -0.999900 | 1.999800 | 0.000006 | 1.684214 |
| | 10 | -0.999900 | 1.999800 | 0.000006 | 1.184225 |
| | 11 | -0.999900 | 1.999800 | 0.000006 | 0.684241 |
| | 12 | -0.999900 | 1.999800 | 0.000006 | 0.184262 |
| | 13 | -0.519400 | 0.522700 | 0.080782 | 0.000000 |
| 9 | 3 | 0 | 0 | 0.062500 | 6.176471 |
| | 4 | 0 | 0 | 0.065891 | 5.412834 |
| | 5 | 0 | 0 | 0.069383 | 4.710455 |
| | 6 | 0 | 0 | 0.072937 | 4.070034 |
| | 7 | 0 | 0 | 0.076512 | 3.491500 |
| | 8 | 0 | 0 | 0.080057 | 2.973928 |
| | 9 | -0.999900 | 1.999800 | 0.000005 | 2.473931 |
| | 10 | -0.999900 | 1.999800 | 0.000005 | 1.973937 |
| | 11 | -0.999900 | 1.999800 | 0.000005 | 1.473947 |
| | 12 | -0.999900 | 1.999800 | 0.000005 | 0.973962 |
| | 13 | -0.999900 | 1.999800 | 0.000005 | 0.473980 |
| | 14 | -0.992200 | 1.879900 | 0.006103 | 0.000000 |
| 10 | 3 | 0 | 0 | 0.055556 | 7.157895 |
| | 4 | 0 | 0 | 0.058282 | 6.369489 |
| | 5 | 0 | 0 | 0.061089 | 5.636078 |
| | 6 | 0 | 0 | 0.063955 | 4.958505 |
| | 7 | 0 | 0 | 0.066852 | 4.337081 |
| | 8 | 0 | 0 | 0.069750 | 3.771515 |
| | 9 | 0 | 0 | 0.072614 | 3.260873 |
| | 10 | -0.999900 | 1.999800 | 0.000004 | 2.760875 |
| | 11 | -0.999900 | 1.999800 | 0.000004 | 2.260881 |
| | 12 | -0.999900 | 1.999800 | 0.000004 | 1.760891 |
| | 13 | -0.999900 | 1.999800 | 0.000005 | 1.260904 |
| | 14 | -0.999900 | 1.999800 | 0.000005 | 0.760921 |
| | 15 | -0.999900 | 1.999800 | 0.000005 | 0.260942 |
| | 16 | -0.944400 | 0.977100 | 0.049090 | 0.000000 |
| 11 | 3 | 0 | 0 | 0.050000 | 8.142857 |
| | 4 | 0 | 0 | 0.052239 | 7.334347 |
| | 5 | 0 | 0 | 0.054542 | 6.575661 |
| | 6 | 0 | 0 | 0.056897 | 5.867674 |
| | 7 | 0 | 0 | 0.059285 | 5.210882 |
| | 8 | 0 | 0 | 0.061687 | 4.605349 |
| | 9 | 0 | 0 | 0.064081 | 4.050667 |
| | 10 | 0 | 0 | 0.066442 | 3.545934 |
| | 11 | -0.999900 | 1.999800 | 0.000004 | 3.045936 |
| | 12 | -0.999900 | 1.999800 | 0.000004 | 2.545942 |
| | 13 | -0.999900 | 1.999800 | 0.000004 | 2.045951 |
| | 14 | -0.999900 | 1.999800 | 0.000004 | 1.545963 |
| | 15 | -0.999900 | 1.999800 | 0.000004 | 1.045979 |
| | 16 | -0.999900 | 1.999800 | 0.000004 | 0.545998 |
| | 17 | -0.999900 | 1.999800 | 0.000004 | 0.046021 |
| | 18 | 0.155500 | 0.070700 | 0.082888 | 0.000000 |
| 12 | 3 | 0 | 0 | 0.045455 | 9.130435 |
| | 4 | 0 | 0 | 0.047325 | 8.305287 |
| | 5 | 0 | 0 | 0.049248 | 7.525642 |
| | 6 | 0 | 0 | 0.051215 | 6.792350 |
| | 7 | 0 | 0 | 0.053213 | 6.105989 |
| | 8 | 0 | 0 | 0.055230 | 5.466819 |
| | 9 | 0 | 0 | 0.057251 | 4.874750 |
| | 10 | 0 | 0 | 0.059260 | 4.329316 |
| | 11 | -0.999900 | 1.999800 | 0.000004 | 3.829316 |
| | 12 | -0.999900 | 1.999800 | 0.000004 | 3.329318 |
| | 13 | -0.999900 | 1.999800 | 0.000004 | 2.829324 |
| | 14 | -0.999900 | 1.999800 | 0.000004 | 2.329332 |
| | 15 | -0.999900 | 1.999800 | 0.000004 | 1.829344 |
| | 16 | -0.999900 | 1.999800 | 0.000004 | 1.329358 |
| | 17 | -0.999900 | 1.999800 | 0.000004 | 0.829376 |
| | 18 | -0.999900 | 1.999800 | 0.000004 | 0.329397 |
| | 19 | -0.970100 | 1.269300 | 0.029236 | 0.000000 |
| 13 | 3 | 0 | 0 | 0.041667 | 10.120000 |
| | 4 | 0 | 0 | 0.043253 | 9.280862 |
| | 5 | 0 | 0 | 0.044882 | 8.483567 |
| | 6 | 0 | 0 | 0.046547 | 7.728921 |
| | 7 | 0 | 0 | 0.048242 | 7.017526 |
| | 8 | 0 | 0 | 0.049956 | 6.349748 |
| | 9 | 0 | 0 | 0.051680 | 5.725686 |
| | 10 | 0 | 0 | 0.053403 | 5.145155 |
| | 11 | 0 | 0 | 0.055111 | 4.607666 |
| | 12 | -0.999900 | 1.999800 | 0.000003 | 4.107667 |
| | 13 | -0.999900 | 1.999800 | 0.000003 | 3.607669 |
| | 14 | -0.999900 | 1.999800 | 0.000003 | 3.107675 |
| | 15 | -0.999900 | 1.999800 | 0.000003 | 2.607683 |
| | 16 | -0.999900 | 1.999800 | 0.000003 | 2.107693 |
| | 17 | -0.999900 | 1.999800 | 0.000004 | 1.607707 |
| | 18 | -0.999900 | 1.999800 | 0.000004 | 1.107724 |
| | 19 | -0.999900 | 1.999800 | 0.000004 | 0.607743 |
| | 20 | -0.999900 | 1.999800 | 0.000004 | 0.107766 |
| | 21 | -0.213800 | 0.247100 | 0.060848 | 0.000000 |
| 14 | 3 | 0 | 0 | 0.038462 | 11.111111 |
| | 4 | 0 | 0 | 0.039823 | 10.260047 |
| | 5 | 0 | 0 | 0.041220 | 9.447694 |
| | 6 | 0 | 0 | 0.042648 | 8.674803 |
| | 7 | 0 | 0 | 0.044102 | 7.941973 |
| | 8 | 0 | 0 | 0.045575 | 7.249620 |
| | 9 | 0 | 0 | 0.047060 | 6.597958 |
| | 10 | 0 | 0 | 0.048549 | 5.986978 |
| | 11 | 0 | 0 | 0.050033 | 5.416436 |
| | 12 | 0 | 0 | 0.051503 | 4.885840 |
| | 13 | -0.999900 | 1.999800 | 0.000003 | 4.385841 |
| | 14 | -0.999900 | 1.999800 | 0.000003 | 3.885844 |
| | 15 | -0.999900 | 1.999800 | 0.000003 | 3.385849 |
| | 16 | -0.999900 | 1.999800 | 0.000003 | 2.885856 |
| | 17 | -0.999900 | 1.999800 | 0.000003 | 2.385867 |
| | 18 | -0.999900 | 1.999800 | 0.000003 | 1.885880 |
| | 19 | -0.999900 | 1.999800 | 0.000003 | 1.385895 |
| | 20 | -0.999900 | 1.999800 | 0.000003 | 0.885914 |
| | 21 | -0.999900 | 1.999800 | 0.000003 | 0.385935 |
| | 22 | -0.966500 | 1.498700 | 0.016757 | 0.000000 |
| 15 | 3 | 0 | 0 | 0.035714 | 12.103448 |
| | 4 | 0 | 0 | 0.036896 | 11.242099 |
| | 5 | 0 | 0 | 0.038107 | 10.416751 |
| | 6 | 0 | 0 | 0.039344 | 9.628104 |
| | 7 | 0 | 0 | 0.040604 | 8.876733 |
| | 8 | 0 | 0 | 0.041882 | 8.163076 |
| | 9 | 0 | 0 | 0.043172 | 7.487410 |
| | 10 | 0 | 0 | 0.044469 | 6.849842 |
| | 11 | 0 | 0 | 0.045767 | 6.250290 |
| | 12 | 0 | 0 | 0.047058 | 5.688473 |
| | 13 | 0 | 0 | 0.048336 | 5.163911 |
| | 14 | -0.999900 | 1.999800 | 0.000003 | 4.663912 |
| | 15 | -0.999900 | 1.999800 | 0.000003 | 4.163915 |
| | 16 | -0.999900 | 1.999800 | 0.000003 | 3.663920 |
| | 17 | -0.999900 | 1.999800 | 0.000003 | 3.163927 |
| | 18 | -0.999900 | 1.999800 | 0.000003 | 2.663937 |
| | 19 | -0.999900 | 1.999800 | 0.000003 | 2.163950 |
| | 20 | -0.999900 | 1.999800 | 0.000003 | 1.663964 |
| | 21 | -0.999900 | 1.999800 | 0.000003 | 1.163982 |
| | 22 | -0.999900 | 1.999800 | 0.000003 | 0.664002 |
| | 23 | -0.999900 | 1.999800 | 0.000003 | 0.164025 |
| | 24 | -0.445100 | 0.455900 | 0.045292 | 0.000000 |
| 16 | 3 | 0 | 0 | 0.033333 | 13.096774 |
| | 4 | 0 | 0 | 0.034368 | 12.226463 |
| | 5 | 0 | 0 | 0.035428 | 11.389793 |
| | 6 | 0 | 0 | 0.036510 | 10.587406 |
| | 7 | 0 | 0 | 0.037612 | 9.819853 |
| | 8 | 0 | 0 | 0.038730 | 9.087570 |
| | 9 | 0 | 0 | 0.039860 | 8.390873 |
| | 10 | 0 | 0 | 0.040999 | 7.729937 |
| | 11 | 0 | 0 | 0.042141 | 7.104789 |
| | 12 | 0 | 0 | 0.043281 | 6.515297 |
| | 13 | 0 | 0 | 0.044414 | 5.961164 |
| | 14 | 0 | 0 | 0.045535 | 5.441925 |
| | 15 | -0.999900 | 1.999800 | 0.000003 | 4.941926 |
| | 16 | -0.999900 | 1.999800 | 0.000003 | 4.441929 |
| | 17 | -0.999900 | 1.999800 | 0.000003 | 3.941934 |
| | 18 | -0.999900 | 1.999800 | 0.000003 | 3.441941 |
| | 19 | -0.999900 | 1.999800 | 0.000003 | 2.941951 |
| | 20 | -0.999900 | 1.999800 | 0.000003 | 2.441962 |
| | 21 | -0.999900 | 1.999800 | 0.000003 | 1.941977 |
| | 22 | -0.999900 | 1.999800 | 0.000003 | 1.441993 |
| | 23 | -0.999900 | 1.999800 | 0.000003 | 0.942012 |
| | 24 | -0.999900 | 1.999800 | 0.000003 | 0.442033 |
| | 25 | -0.984500 | 1.744600 | 0.007452 | 0.000000 |
| | 3 | 0 | 0 | 0.031250 | 14.090909 |
| | 4 | 0 | 0 | 0.032164 | 13.212722 |
| | 5 | 0 | 0 | 0.033099 | 12.366098 |
| | 6 | 0 | 0 | 0.034053 | 11.551631 |
| | 7 | 0 | 0 | 0.035024 | 10.769837 |
| | 8 | 0 | 0 | 0.036010 | 10.021148 |
| | 9 | 0 | 0 | 0.037008 | 9.305892 |
| | 10 | 0 | 0 | 0.038014 | 8.624289 |
| | 11 | 0 | 0 | 0.039025 | 7.976438 |
| | 12 | 0 | 0 | 0.040038 | 7.362309 |
| | 13 | 0 | 0 | 0.041047 | 6.781738 |
| | 14 | 0 | 0 | 0.042049 | 6.234421 |
| 17 | 15 | 0 | 0 | 0.043040 | 5.719911 |
| | 16 | -0.999900 | 1.999800 | 0.000003 | 5.219912 |
| | 17 | -0.999900 | 1.999800 | 0.000003 | 4.719915 |
| | 18 | -0.999900 | 1.999800 | 0.000003 | 4.219920 |
| | 19 | -0.999900 | 1.999800 | 0.000003 | 3.719927 |
| | 20 | -0.999900 | 1.999800 | 0.000003 | 3.219936 |
| | 21 | -0.999900 | 1.999800 | 0.000003 | 2.719947 |
| | 22 | -0.999900 | 1.999800 | 0.000003 | 2.219961 |
| | 23 | -0.999900 | 1.999800 | 0.000003 | 1.719976 |
| | 24 | -0.999900 | 1.999800 | 0.000003 | 1.219994 |
| | 25 | -0.999900 | 1.999800 | 0.000003 | 0.720015 |
| | 26 | -0.999900 | 1.999800 | 0.000003 | 0.220037 |
| | 27 | -0.743900 | 0.747300 | 0.033369 | 0.000000 |
| | 3 | 0 | 0 | 0.029412 | 15.085714 |
| | 4 | 0 | 0 | 0.030225 | 14.200551 |
| | 5 | 0 | 0 | 0.031055 | 13.345110 |
| | 6 | 0 | 0 | 0.031903 | 12.519940 |
| | 7 | 0 | 0 | 0.032765 | 11.725526 |
| | 8 | 0 | 0 | 0.033641 | 10.962283 |
| 18 | 9 | 0 | 0 | 0.034528 | 10.230544 |
| | 10 | 0 | 0 | 0.035423 | 9.530554 |
| | 11 | 0 | 0 | 0.036323 | 8.862459 |
| | 12 | 0 | 0 | 0.037227 | 8.226302 |
| | 13 | 0 | 0 | 0.038130 | 7.622011 |
| | 14 | 0 | 0 | 0.039029 | 7.049401 |
| | 15 | 0 | 0 | 0.039921 | 6.508168 |
| | 16 | 0 | 0 | 0.040803 | 5.997887 |
| | 17 | -0.999900 | 1.999800 | 0.000002 | 5.497888 |
| | 18 | -0.999900 | 1.999800 | 0.000002 | 4.997891 |
| | 19 | -0.999900 | 1.999800 | 0.000002 | 4.497896 |
| | 20 | -0.999900 | 1.999800 | 0.000002 | 3.997903 |
| | 21 | -0.999900 | 1.999800 | 0.000003 | 3.497911 |
| | 22 | -0.999900 | 1.999800 | 0.000003 | 2.997922 |
| | 23 | -0.999900 | 1.999800 | 0.000003 | 2.497935 |
| | 24 | -0.999900 | 1.999800 | 0.000003 | 1.997950 |
| | 25 | -0.999900 | 1.999800 | 0.000003 | 1.497967 |
| | 26 | -0.999900 | 1.999800 | 0.000003 | 0.997986 |
| | 27 | -0.999900 | 1.999800 | 0.000003 | 0.498008 |
| | 28 | -0.998100 | 1.989800 | 0.000228 | 0.000000 |
| 19 | 3 | 0 | 0 | 0.027778 | 16.081081 |
| | 4 | 0 | 0 | 0.028505 | 15.189695 |
| | 5 | 0 | 0 | 0.029249 | 14.326392 |
| | 6 | 0 | 0 | 0.030006 | 13.491676 |
| | 7 | 0 | 0 | 0.030777 | 12.686002 |
| | 8 | 0 | 0 | 0.031560 | 11.909768 |
| | 9 | 0 | 0 | 0.032352 | 11.163303 |
| | 10 | 0 | 0 | 0.033153 | 10.446864 |
| | 11 | 0 | 0 | 0.033959 | 9.760628 |
| | 12 | 0 | 0 | 0.034770 | 9.104685 |
| | 13 | 0 | 0 | 0.035581 | 8.479035 |
| | 14 | 0 | 0 | 0.036391 | 7.883577 |
| | 15 | 0 | 0 | 0.037197 | 7.318113 |
| | 16 | 0 | 0 | 0.037997 | 6.782343 |
| | 17 | 0 | 0 | 0.038786 | 6.275864 |
| | 18 | -0.999900 | 1.999800 | 0.000002 | 5.775865 |
| | 19 | -0.999900 | 1.999800 | 0.000002 | 5.275868 |
| | 20 | -0.999900 | 1.999800 | 0.000002 | 4.775873 |
| | 21 | -0.999900 | 1.999800 | 0.000002 | 4.275880 |
| | 22 | -0.999900 | 1.999800 | 0.000002 | 3.775888 |
| | 23 | -0.999900 | 1.999800 | 0.000002 | 3.275899 |
| | 24 | -0.999900 | 1.999800 | 0.000002 | 2.775911 |
| | 25 | -0.999900 | 1.999800 | 0.000002 | 2.275925 |
| | 26 | -0.999900 | 1.999800 | 0.000002 | 1.775942 |
| | 27 | -0.999900 | 1.999800 | 0.000003 | 1.275960 |
| | 28 | -0.999900 | 1.999800 | 0.000003 | 0.775981 |
| | 29 | -0.999900 | 1.999800 | 0.000003 | 0.276004 |
| | 30 | -0.888700 | 1.023500 | 0.023926 | 0.000000 |
| 20 | 3 | 0 | 0 | 0.026316 | 17.076923 |
| | 4 | 0 | 0 | 0.026971 | 16.179953 |
| | 5 | 0 | 0 | 0.027640 | 15.309594 |
| | 6 | 0 | 0 | 0.028321 | 14.466313 |
| | 7 | 0 | 0 | 0.029014 | 13.650536 |
| | 8 | 0 | 0 | 0.029717 | 12.862638 |
| | 9 | 0 | 0 | 0.030430 | 12.102943 |
| | 10 | 0 | 0 | 0.031150 | 11.371713 |
| | 11 | 0 | 0 | 0.031876 | 10.669142 |
| | 12 | 0 | 0 | 0.032606 | 9.995354 |
| | 13 | 0 | 0 | 0.033338 | 9.350396 |
| | 14 | 0 | 0 | 0.034071 | 8.734233 |
| | 15 | 0 | 0 | 0.034802 | 8.146746 |
| | 16 | 0 | 0 | 0.035528 | 7.587730 |
| | 17 | 0 | 0 | 0.036248 | 7.056892 |
| | 18 | 0 | 0 | 0.036959 | 6.553850 |
| | 19 | -0.999900 | 1.999800 | 0.000002 | 6.053852 |
| | 20 | -0.999900 | 1.999800 | 0.000002 | 5.553855 |
| | 21 | -0.999900 | 1.999800 | 0.000002 | 5.053860 |
| | 22 | -0.999900 | 1.999800 | 0.000002 | 4.553866 |
| | 23 | -0.999900 | 1.999800 | 0.000002 | 4.053874 |
| | 24 | -0.999900 | 1.999800 | 0.000002 | 3.553884 |
| | 25 | -0.999900 | 1.999800 | 0.000002 | 3.053896 |
| | 26 | -0.999900 | 1.999800 | 0.000002 | 2.553910 |
| | 27 | -0.999900 | 1.999800 | 0.000002 | 2.053926 |
| | 28 | -0.999900 | 1.999800 | 0.000002 | 1.553943 |
| | 29 | -0.999900 | 1.999800 | 0.000002 | 1.053963 |
| | 30 | -0.999900 | 1.999800 | 0.000002 | 0.553985 |
| | 31 | -0.999900 | 1.999800 | 0.000002 | 0.054008 |
| | 32 | -0.102400 | 0.114000 | 0.044720 | 0.000000 |